\documentclass{amsart}
\usepackage[all]{xy}

\theoremstyle{plain}
\newtheorem{theorem}{Theorem}[section]
\newtheorem{proposition}[theorem]{Proposition}
\newtheorem{lemma}[theorem]{Lemma}
\newtheorem{corollary}[theorem]{Corollary}
\newtheorem{conjecture}[theorem]{Conjecture}

\theoremstyle{definition}
\newtheorem{definition}[theorem]{Definition}
\newtheorem{example}[theorem]{Example}

\theoremstyle{remark}
\newtheorem{remark}[theorem]{Remark}
\newtheorem{remarks}[theorem]{Remarks}

\newcommand{\secref}[1]{Section~\ref{#1}}
\newcommand{\thmref}[1]{Theorem~\ref{#1}}
\newcommand{\propref}[1]{Proposition~\ref{#1}}
\newcommand{\lemref}[1]{Lemma~\ref{#1}}
\newcommand{\corref}[1]{Corollary~\ref{#1}}
\newcommand{\conjref}[1]{Conjecture~\ref{#1}}

\newcommand{\exref}[1]{Example~\ref{#1}}

\newcommand{\remref}[1]{Remark~\ref{#1}}

\newcommand{\defref}[1]{Definition~\ref{#1}}

\def\R{{\mathbb \R}}

\def\Q{{\mathbb Q}}

\def\H{{\mathbb H}}

\def\map{\mathrm{map}}
\def\Der{\mathrm{Der}}
\def\Hom{\mathrm{Hom}}
\def\Rel{\mathrm{Rel}}
\def\ker{\mathrm{ker}}
\def\im{\mathrm{im}}

\begin{document}

\title[Rationalized Evaluation Subgroups]{Rationalized Evaluation Subgroups of a Map and the Rationalized $G$-Sequence}

\author{Gregory  Lupton}

\address{Department of Mathematics,
          Cleveland State University,
          Cleveland OH 44115}

\email{G.Lupton@csuohio.edu}

\author{Samuel Bruce Smith}

\address{Department of Mathematics,
  Saint Joseph's University,
  Philadelphia PA 19131}

\email{smith@sju.edu}

\date{\today}



\begin{abstract}
Let $f\colon X \rightarrow Y$ be a based map of simply connected
spaces.   The corresponding evaluation map $\omega \colon
\map(X,Y;f) \to Y$ induces a homomorphism of homotopy groups whose
image in $\pi_n(Y)$ is called the $n$th evaluation subgroup of
$f$. The $n$th Gottlieb group of $X$ occurs as the special case in
which $Y = X$ and $f = 1_X$.  We identify the homomorphism induced
on rational homotopy groups by this evaluation map, in terms of a
map of complexes of derivations constructed using Sullivan minimal
models.  Our identification allows for the characterization of the
rationalization of the $n$th evaluation subgroup of $f$. It also
allows for the identification of several long exact sequences of
rational homotopy groups, including the long exact sequence
induced on rational homotopy groups by the evaluation fibration.
As a consequence, we obtain an identification of the
rationalization of the so-called $G$-sequence of the map $f$. This
is a sequence---in general not exact---of groups and homomorphisms
that includes the Gottlieb groups of $X$ and the evaluation
subgroups of $f$. We use these results to study the $G$-sequence
in the context of rational homotopy theory.     We give new
examples of non-exact $G$-sequences and  uncover a relationship
between the homology of the rational $G$-sequence and negative
derivations of rational cohomology. We also  relate the splitting
of the rational $G$-sequence of a fibre inclusion to a well-known
conjecture in rational homotopy theory.
\end{abstract}

\maketitle

\section{Introduction}\label{sec:intro}

Suppose given a based map $f \colon X \rightarrow Y$ of simply
connected CW complexes.  Denote by $\map(X, Y; f)$ the path
component of the space of (unbased) maps $X\to Y$ consisting of
those maps that are homotopic to $f$. Then evaluation at the
basepoint of $X$ gives a based map $\omega \colon \map(X, Y; f)
\rightarrow Y$. We refer to this map as the \emph{evaluation map}.
We define the \emph{$n$th evaluation subgroup of $f$} to be the
subgroup $G_{n}(Y, X; f) = \omega_{\#}(\pi_{n}(\map(X,Y;f)))$ of
$\pi_n(Y)$.  The famous Gottlieb groups $G_{*}(X)$ occur as the
special case in which $X=Y$ and $f= 1_{X}$ \cite{Go1}.  The
Gottlieb groups of a space have been much studied by homotopy
theorists (see \cite{Oprea} for a survey of results and
references). While many general results are known, explicit
computation of $G_*(X)$ appears difficult and is limited to a
small number of sporadic examples. One reason that accounts in
part for this difficulty is the fact that a map of spaces $f
\colon X \to Y$ does not necessarily induce a corresponding
homomorphism of Gottlieb groups, since in general $f_\#(G_n(X))
\not\subseteq G_n(Y)$.  In particular, attempts to study $G_*(X)$
via a cell decomposition of $X$ are frustrated, since it is not
clear what effect a cell attachment may have on the Gottlieb
groups.  One tool for studying $G_*(X)$ that attempts to
circumvent this problem is the so-called \emph{$G$-sequence of a
map}, introduced by Woo and Lee in \cite{L-W1}.

The $G$-sequence of a map $f \colon X \rightarrow Y$ is a sequence
$$\cdots \rightarrow G_{n}(X) \rightarrow G_{n}(Y, X;f) \rightarrow
G_{n}^{rel}(Y,X;f) \rightarrow G_{n-1}(X) \rightarrow \cdots $$
of groups and homomorphisms that derives from the long exact
homotopy sequence of the map $f$. In this sequence, $G_n(X)$ is
the Gottlieb group of $X$, $G_{n}(Y, X;f)$ is the evaluation
subgroup of $f$ introduced above, and the third term
$G_{n}^{rel}(Y,X;f)$ is a suitably defined  ``relative" term. The
$G$-sequence arises as follows. One has the following commutative
diagram of spaces:
\begin{displaymath}
\xymatrix{\map(X, X; 1) \ar[r]^{f_*} \ar[d]_{\omega}& \map(X, Y;
f)\ar[d]^{\omega}\\
X \ar[r]_{f} & Y}
\end{displaymath}
Now pass to the corresponding induced homomorphisms of homotopy
groups.  In a standard way, the induced homomorphisms of homotopy
groups $(f_*)_\#$ and $f_\#$ can be fitted into the long exact
homotopy sequences of the maps  $f_*$ and $f$ respectively.  Then
the evaluation maps induce maps of each term in the long exact
sequences, resulting in a homotopy ladder. The $G$-sequence of the
map $f$ is then the image of the top long exact homotopy sequence
in that of the bottom.  A portion of the $G$-sequence is shown
here:
$$
\xymatrix{
  \pi_{n+1}(f_*) \ar[r]
  \ar@{->>}[d]^{\omega_{\#}} & \pi_{n}(\map(X,
X;1)) \ar@{->>}[d]^{\omega_{\#}}
  \ar[r]^{(f_*)_{\#}}
  & \pi_{n}(\map(X, Y;f))   \ar@{->>}[d]^{\omega_{\#}}
    \\
  G^{rel}_{n+1}(Y, X;f) \ar@{^{(}->}[d] \ar[r] &
   G_{n}(X) \ar@{^{(}->}[d] \ar[r]^{f_{\#}}
  & G_{n}(Y, X;f) \ar@{^{(}->}[d]
     \\
     \pi_{n+1}(f) \ar[r] & \pi_{n}(X)  \ar[r]^{f_{\#}}  & \pi_{n}(Y)
      }
$$
The homomorphisms in the $G$-sequence are restrictions of those in
the long exact homotopy sequence of the map $f$. Therefore, the
$G$-sequence forms a chain complex (consecutive compositions are
trivial).  Some general conditions are known under which the
$G$-sequence is exact (e.g.~\cite{L-W4, P-W}), but in general it
is not exact (e.g.~\cite{L-W4}).

In this paper we bring the techniques of rational homotopy theory
to bear on problems and questions concerning evaluation subgroups
(of a map) in general, and the $G$-sequence in particular.  Our
main goal is to expand the range of application of these
techniques in this area.  To this end, we are primarily concerned
with establishing a suitable framework for considering such
questions. At the same time, we obtain a number of results of
interest in their own right.

The paper is organized as follows:  Our main results are
established in Sections $2$ and $3$. The basic result in
\secref{sec:deriv spaces} is \thmref{thm:equivalence of squares},
in which we identify the map induced on rational homotopy groups
by the evaluation map $\omega \colon \map(X, Y; f) \rightarrow Y$.
We describe this induced homomorphism as the map induced on
homology by a map of complexes of derivations of the Sullivan
minimal models of $X$ and $Y$. \thmref{thm:equivalence of squares}
has a number of immediate corollaries.  For instance, we obtain a
characterization of the rationalized evaluation subgroups of a map
(\corref{cor:rational eval subgp of f}) that extends a well known
characterization, due to F{\'e}lix and Halperin, of the
(rationalized) Gottlieb groups of a space in terms of derivations
of its minimal model.

In \secref{sec:long exact sequences}, we extend and amplify the
basic result of \thmref{thm:equivalence of squares}.   We show
that several related long exact sequences of rational homotopy
groups are naturally expressed as long exact sequences in homology
of derivation complexes of Sullivan minimal models. In particular,
we obtain a description, within the framework of derivation
spaces, of the $G$-sequence of a map after rationalization.  The
results of this section are of interest in their own right and in
several cases (e.g.~\thmref{thm:long exact seq of eval fib}) they
are independent of the $G$-sequence.

In \secref{sec:Examples}, we use the framework established in
Sections $2$ and $3$ to study questions concerning the
rationalized $G$-sequence. Just as in the integral setting, the
rationalized $G$-sequence is not exact in general.
\exref{ex:non-exact G-seq} gives a simple example for which the
rationalized $G$-sequence fails to be exact at each of the three
types of term that occur. Since non-exactness rationally implies
non-exactness integrally, this example provides a new, complete
example of the failure of exactness of the $G$-sequence. This
example also shows that the framework established in Sections $2$
and $3$ provides an effective setting in which to carry out
explicit computations.  By way of contrast, in \thmref{thm:exact
G-seq} we give one set of conditions under which the rationalized
$G$-sequence is exact.  In \thmref{thm:H^omega=0} we show that
under certain circumstances, the rationalized $G$-sequence may be
exact at all occurrences of one type of term, while failing to be
exact at the other types of term. This same result establishes a
relationship, albeit under rather restrictive circumstances,
between the (vanishing of a certain type of the) $\omega$-homology
of a map $f \colon X \to Y$ and the (vanishing of) negative-degree
derivations of the rational cohomology of $X$. The last
development that we present in \secref{sec:Examples} is a
connection between the $G$-sequence of certain fibre inclusion
maps and a well-known conjecture in rational homotopy theory.  Let
$X \to E \to S^{2r+1}$ be a fibration with base an odd-dimensional
sphere. For certain types of fibre space $X$, a conjecture of
Halperin asserts that the fibration should be rationally TNCZ (see
\conjref{conj:Halperin}
 below). In \thmref{thm:Halperin G-trivial}, we show that
this is the case exactly when the rationalized $G$-sequence of the
fibre inclusion $X \to E$ reduces to a certain short exact
sequence. In this way, we obtain an equivalent phrasing of the
conjecture of Halperin, in terms of the ideas studied in this
paper.

Finally, in a technical appendix, we give full details for several
results from rational homotopy theory used in the proof of
\thmref{thm:equivalence of squares}.

We next discuss existing results on rational homotopy, function
spaces, and the Gottlieb groups of a space that relate to this
paper. Our aim in this discussion is to indicate the basic results
that exist in the area.  At the same time, we identify how these
results relate to our work, and our points of departure from them.
From results in rational homotopy theory, we can identify three
``tributary streams" that flow into this work precisely at
\thmref{thm:equivalence of squares}---our first main result---at
varying levels of generality.  There is also a fourth tributary
flowing into our work, from outside rational homotopy theory, that
includes the work on the $G$-sequence.

The first tributary of antecedent results in rational homotopy
theory concerns the rationalized Gottlieb groups of a space. In
\cite{F-H}, F{\'e}lix and Halperin gave a characterization of the
rationalized Gottlieb groups of a space, in terms of derivations
of the Sullivan minimal model.  At its most specialized level,
\thmref{thm:equivalence of squares} retrieves this
characterization (\corref{cor:rational Gottlieb group}) and
extends it to a similar characterization of the rationalized
evaluation subgroups of a general map (\corref{cor:rational eval
subgp of f}).  F{\'e}lix and Halperin went on to prove a
remarkable result concerning the rationalized Gottlieb groups of a
finite complex \cite[Th.III]{F-H}.   Their result significantly
extends results of Gottlieb from \cite{Go1}, and relates the
rationalized Gottlieb groups with the rational
Lusternik-Schnirelmann category. Unfortunately, no analogous
result seems forthcoming for the rationalized evaluation subgroups
of a map.  Nonetheless, our characterization of the rationalized
evaluation subgroups of a map is as effective for concrete
computations as is the earlier characterization of the
rationalized Gottlieb groups.

The second tributary concerns the rational homotopy type of
$B\,\mathrm{aut}_1(X)$---the classifying space for fibrations with
fibre $X$.   Although we are not concerned with this classifying
space as such, the connection arises because we have an
isomorphism of homotopy groups
$\pi_{i+1}\big(B\,\mathrm{aut}_1(X)\big) \cong \pi_i\big(\map(X,
X; 1)\big)$, and also the Gottlieb groups of $X$ are obtained as
the image in homotopy groups of the connecting homomorphism of the
corresponding classifying fibration \cite[Th.2.6]{Go1}.  At a
higher level of generality than the rationalized Gottlieb groups,
\thmref{thm:equivalence of squares} includes an identification of
the map induced on rational homotopy groups by the evaluation map
$\omega \colon \map(X, X; 1) \rightarrow X$, and in particular it
identifies the rational homotopy groups of $\map(X, X; 1)$. Whilst
the identification of these groups we give is familiar in rational
homotopy theory, our proof for this---obtained by restricting that
of \thmref{thm:equivalence of squares}---is the first direct and
detailed one that has actually appeared in the literature.  We
support this assertion as follows: In \cite[Sec.11]{Su}, Sullivan
sketched (with no proof) a model for the rational homotopy type of
$B\,\mathrm{aut}_1(X)$, from which the description of the rational
homotopy groups of $\map(X, X; 1)$ contained in
\thmref{thm:equivalence of squares} follows. A justification of
Sullivan's model for $B\,\mathrm{aut}_1(X)$ may be gleaned by
collating results from a number of sources spread through the
literature (e.g.~\cite{Sch-St, Tan84, Tan83, Gats}). But to date,
no direct proof of Sullivan's model has been given. Even amongst
those articles that focus specifically on the rational homotopy
\emph{groups}, either of $B\,\mathrm{aut}_1(X)$ or of $\map(X, X;
1)$---and that therefore avoid the technical problems of dealing
with the rational homotopy \emph{type}---we still do not find
complete details. Meier \cite[(1.4), (2.6)]{Me} outlines the basic
idea, but is actually focussed on a special kind of situation in
which the minimal model can be replaced by its cohomology. Grivel
\cite{Griv} focusses on the same special case as Meier, and quotes
Sullivan's model directly. F\'{e}lix and Thomas
\cite[Sec.2.3]{F-T} give exactly the description of the rational
homotopy groups of $\map(X, X; 1)$ contained in
\thmref{thm:equivalence of squares}, but no details of the proof
are given.  As well as including a direct and detailed proof for
the case of $\map(X, X; 1)$, \thmref{thm:equivalence of squares}
extends this identification of the rational homotopy groups to the
general case of $\map(X, Y; f)$.   It thus provides a natural
framework for the study of rational homotopy groups of function
spaces.

Finally, the third tributary from rational homotopy theory
consists of a model, due to Sullivan and Haefliger, for the
rational homotopy type of $\map(X, Y; f)$.  At its full level of
generality, \thmref{thm:equivalence of squares} identifies the map
induced on rational homotopy groups by $\omega \colon \map(X, Y;
f) \rightarrow Y$, a general evaluation map.  The precursor to the
identification we give, and to our general line of proof, is the
approach of Thom in \cite{Th}, although obviously there is no
reference to minimal models in his work. Pursuing the river
analogy a little further, we may think of Thom's result, if not as
the source, then at least as somewhere in the headwaters. We
emphasize the connection between \thmref{thm:equivalence of
squares} and Thom's approach by retrieving a basic result of his
in \corref{cor:Thom}. Coming further downstream, Sullivan also
described in \cite[Sec.11]{Su} a model for the space of sections
of a fibration homotopic to a given section. By specializing to
the trivial fibration $X \times Y \to X$, this yields a model for
the function space $\map(X, Y; f)$---and more generally a model
for the rational homotopy \emph{type} of the general evaluation
map $\omega$. A detailed proof for Sullivan's model in this case
was given by Haefliger \cite{Hae}. Now this model should in
principle determine the rational homotopy groups of the function
space (see \cite{F-T} and \cite{Mo-Ra}, where it is used quite
effectively). However, the model in question is a (non-minimal) DG
algebra model. Therefore, the homomorphism induced by $\omega$ on
rational homotopy groups---which is exactly the information we
require to proceed with our development---is available only
indirectly, at best.   By focussing on the rational homotopy
\emph{groups}---as opposed to the rational homotopy type, we have
arrived in \thmref{thm:equivalence of squares} at an entirely new
characterization of the map induced on rational homotopy groups by
$\omega \colon \map(X, Y; f) \rightarrow Y$.  Furthermore, we have
been able to give a direct proof that avoids many of the technical
complexities of Haefliger's work and is completely independent of
it.

The remaining tributary flowing into our work comes from outside
rational homotopy theory, and concerns classical results on the
Gottlieb groups of a space, and more recent results on evaluation
subgroups of a map and the $G$-sequence.  We have already
mentioned some of the results in this area.  In his original work
on evaluation subgroups, Gottlieb observed that a map of spaces
does not necessarily induce a map of Gottlieb groups, and gave
conditions under which it does \cite[Sec.1]{Go1}. Gottlieb briefly
mentions the evaluation subgroups of a map in \cite{Go1}, but did
not study them as such. A number of basic properties of the
evaluation subgroup of a map are established in \cite{Woo-Kim,
Woo-Kim2, L-W2}. As we mentioned earlier, Woo and Lee introduced
the $G$-sequence of a map in \cite{L-W1}.  Results on the
$G$-sequence basically fall into one of three areas: conditions
under which the $G$-sequence is exact (e.g.~\cite[Th.12]{L-W1} and
\cite{L-W3}), examples of non-exactness (e.g.~\cite{L-W4} and
\cite{Pak-W}), and extensions and generalizations of evaluation
subgroups and the $G$-sequence (e.g.~\cite{L-W5, L-W6}). These
results give the stepping-off point for our work in Sections 3 and
4.  Using our description of the rationalized $G$-sequence
(\thmref{thm:G-seq of f}), we extend the known exactness results
to several new cases. We also give new instances of non-exactness.
More significantly, by focussing on the rational setting, our
methods make the production of such examples straightforward. On
the other hand, our results and examples are not restricted to
amplifying previous results in this area. In
\corref{cor:G-sequence Halperin}, we suggest a different kind of
result that relates properties of the $G$-sequence to the
triviality of a fibration.  Furthermore, some results in Section
3---including \thmref{thm:long exact seq of f_*} and
\thmref{thm:long exact seq of eval fib}---are of interest
independently of any relation to the $G$-sequence.

We finish this introduction by setting some notation and
terminology. Throughout this paper, $X$ and $Y$ will denote simply
connected CW complexes of finite type. By \emph{vector space} we
mean a rational graded vector space. By \emph{algebra}, we mean
the kind of commutative graded algebras over the rationals that
arise in rational homotopy.  That is, they are non-negatively
graded, connected ($H^0 = \Q$) and usually simply connected ($H^1
= 0$), with cohomology of finite type.  For a vector space $V$, we
denote the free commutative graded algebra generated by $V$ by
$\Lambda V$. We use the acronym DG to denote differential graded:
Thus, DG vector space, DG algebra, and so-forth.  For a DG
algebra, the differential is of degree $+1$. In other situations,
however, particularly when we consider the complex of derivations
of a DG algebra, the differential is of degree $-1$. We will
generally refer to a DG vector space whose differential is of
degree $-1$ as a chain complex.  If $f \colon A \to B$ is a map,
either topological or algebraic, then $f^*$ denotes
pre-composition by $f$ and $f_*$ denotes post-composition by $f$.
In any setting in which it is appropriate, we use $H(f)$ to denote
the map induced on homology (or cohomology) by $f$, and $f_\#$ to
denote the map induced on homotopy groups by the map of spaces
$f$. A map of DG algebras is called a \emph{quasi-isomorphism} if
it induces an isomorphism on cohomology.  We use $\omega$ in a
generic way to denote an evaluation map, and we denote the
identity map of a topological space or the identity homomorphism
of an algebra by $1$.  We denote the rationalization of a space
$X$ by $X_\Q$ and of a map $f$ by $f_\Q$ (cf. \cite{H-M-R}).

We assume that the reader is familiar with the basics of rational
homotopy. Our general  reference for this material is
\cite{F-H-T}. We recall, in particular, that a space $X$ has a
\emph{minimal model} $\mathcal{M}_X$, which is a certain type of
DG algebra. Namely, $\mathcal{M}_X$ is a free algebra $\Lambda V$
with a decomposable differential, that is, $d(V) \subseteq
\Lambda^{\geq2} V$. Furthermore, a map of spaces $f \colon X \to
Y$ induces a map of minimal models $\mathcal{M}_{f} \colon
\mathcal{M}_Y \to \mathcal{M}_X$.  We refer to this induced map as
the \emph{Sullivan minimal model} of the map $f$. It is a complete
rational homotopy invariant for a map, and in principle all
rational homotopy theoretic information about $f$ can be retrieved
from it. Passing to cohomology, for example, gives
$H(\mathcal{M}_{f}) \colon H^*(\mathcal{M}_Y) \to
H^*(\mathcal{M}_X)$, which corresponds to the homomorphism of
rational cohomology algebras induced by $f$. The results of this
paper illustrate how deeper information about a space or map may
be retrieved from the minimal model by making correspondingly more
sophisticated constructions with the model.

\section{Derivation Spaces}\label{sec:deriv
spaces}

Our purpose in this and the next section is to give a unified
description in rational homotopy theory  of all  the terms
involved in the definition of  the $G$-sequence.  Informally
stated, we show that the homology theory of derivation complexes
of Sullivan minimal models provides  an algebraic model for the
rational homotopy theory of function spaces at the level of
homotopy groups.

We focus on the following commutative square that appears in the
homotopy ladder from which the $G$-sequence arises:
\begin{equation}\label{eq:homotopy square}
\xymatrix{\pi_n\big(\map(X, X; 1)\big) \ar[r]^{(f_*)_\#}
\ar[d]_{\omega_\#}& \pi_n\big(\map(X, Y;
f)\big)\ar[d]^{\omega_\#}\\
\pi_n(X) \ar[r]_{f_\#} & \pi_n(Y)}
\end{equation}
It turns out that identifying the rationalization of this
commutative square is sufficient to arrive not only at a
characterization of the rationalized evaluation subgroups of $f$,
but also at a description of the rationalization of the
$G$-sequence. Furthermore, our identification of the
rationalization of this square allows us to conclude several
subsidiary results of interest.

We say two maps of vector spaces $f\colon U \to V$ and $g\colon U'
\to V'$ are \emph{equivalent} if there exist isomorphisms $\alpha$
and $\beta$ which make the diagram
\begin{displaymath}
\xymatrix{U \ar[r]^{f} \ar[d]_{\alpha}^{\cong} & V
\ar[d]_{\beta}^{\cong}\\
U' \ar[r]_{g} & V'}
\end{displaymath}
commutative.  This notion of equivalence for vector space maps
extends in the obvious way to sequences of vector space maps,
commutative squares of vector space maps, and any other diagram of
vector space maps.

Suppose a DG algebra $(A, d_A)$ is isomorphic to $\Q$ in degree
zero, that is, $A^0 \cong \Q$.  Then the map $\varepsilon \colon A
\rightarrow \Q$ that sends all elements of positive degree to
zero, and is the identity in degree zero, is an augmentation. This
will be the situation in all cases of interest to us here, and
thus we refer to $\varepsilon \colon A \rightarrow \Q$ as
\emph{the} augmentation.
 Here, as in the sequel, we regard $\Q$ as the
trivial DG algebra concentrated in degree zero and with trivial
differential. Thus $\varepsilon$ is a DG algebra map.

Given DG algebras $(A,d_A)$ and $(B,d_B)$ and a (fixed) DG algebra
map $\phi \colon A \rightarrow B$, define a
\emph{$\phi$-derivation} of degree $n$ to be a linear map $\theta
\colon A \rightarrow B$ that \emph{reduces degree by $n$} and
satisfies the derivation law $\theta (xy) = \theta(x)\phi(y) -
(-1)^{n|x|} \phi(x) \theta(y)$.  We will only consider derivations
of positive degree, that is, those that reduce degree by some
positive integer. When $n=1$ we require additionally that $d_{B}
\circ \theta = - \theta \circ d_{A}$.  Let $\Der_{n}(A, B;
{\phi})$ denote the vector space of $\phi$-derivations of degree
$n$, for $n>0$. Finally, define a linear map $\delta \colon
\Der_{n}(A, B; {\phi}) \to \Der_{n-1}(A, B; {\phi})$ by
$\delta(\theta) = d_B \circ \theta  - (-1)^{|\theta|} \theta \circ
d_A$.   A standard check now shows that $\delta^2 = 0$ and thus
$\big(\Der_{*}(A, B; {\phi}), \delta\big)$ is a chain complex. In
order to cut down on cumbersome notation, we will usually suppress
the differential from our notation, and write $H_n\big(\Der(A, B;
{\phi})\big)$ for the homology in degree $n$ of the chain complex
$\big(\Der_{*}(A, B; {\phi}), \delta\big)$. That is,
$H_n\big(\Der(A, B; {\phi})\big)$ denotes the homology represented
by $\delta$-cycles of $\Der_{*}(A, B; {\phi})$ that reduce degree
by $n$.

A special case of the preceding that is of interest to us is the
one in which $A=B$ and $\phi = 1_B$.  In this case, the chain
complex of derivations $\Der_{*}(B, B; 1)$ is just the usual
complex of derivations on the DG algebra $B$.  Note once again
that we restrict the derivations in degree $1$ to the cycles and
that the complex is zero in non-positive degrees.  Pre-composition
with the DG algebra map $\phi \colon A \to B$ thus gives a map of
chain complexes $\phi^* \colon \Der(B, B; 1) \rightarrow \Der(A,
B; \phi)$. Furthermore, post-composition by the augmentation
$\varepsilon \colon B \to \Q$ induces DG vector space maps
$\varepsilon_* \colon \Der_{*}(A, B; \phi) \to \Der_{*}(A, \Q;
\varepsilon)$ and $\varepsilon_* \colon \Der_{*}(B, B; 1) \to
\Der_{*}(B, \Q; \varepsilon)$.

All of the above can be applied to a map of minimal models.
Suppose $f \colon X \to Y$ is a map of spaces, and
$\mathcal{M}_{f} \colon \mathcal{M}_{Y} \to \mathcal{M}_{X}$ is
the corresponding Sullivan minimal model of the map $f$. Then we
have a commutative square of chain complexes
\begin{equation}\label{eq:derivation square}
\xymatrix{ \Der_*(\mathcal{M}_{X}, \mathcal{M}_{X}; 1)
\ar[r]^{(\mathcal{M}_{f})^*} \ar[d]_{\varepsilon_*} &
\Der_*(\mathcal{M}_{Y}, \mathcal{M}_{X}; \mathcal{M}_{f})
\ar[d]^{\varepsilon_*}\\
\Der_*(\mathcal{M}_{X}, \Q; \varepsilon)
\ar[r]^{\widehat{(\mathcal{M}_{f})^*}} & \Der_*(\mathcal{M}_{Y},
\Q; \varepsilon)}
\end{equation}
In this square, $\varepsilon \colon \mathcal{M}_{X} \to \Q$ is the
augmentation.    Both horizontal maps are obtained by
pre-composing with the same map $\mathcal{M}_{f}$, but in
different contexts.  Since we will need to distinguish between
these two maps notationally in the sequel, we have used an extra
decoration on the bottom one.

The main result of this section is the following:

\begin{theorem}\label{thm:equivalence of squares}
Let $X$ and $Y$ be simply connected CW complexes of finite type,
with $X$ finite.  For $n \geq 2$, the commutative square obtained
by rationalizing (\ref{eq:homotopy square}) is equivalent to the
square obtained by passing to homology in degree $n$ from
(\ref{eq:derivation square}). That is, the commutative squares
\begin{displaymath}
\xymatrix{\pi_n\big(\map(X, X; 1)\big)\otimes\Q
\ar[rr]^{(f_*)_\#\otimes\Q} \ar[d]_{\omega_\#\otimes\Q} & &
\pi_n\big(\map(X, Y;
f)\big)\otimes\Q\ar[d]^{\omega_\#\otimes\Q}\\
\pi_n(X)\otimes\Q \ar[rr]_{f_\#\otimes\Q} & & \pi_n(Y)\otimes\Q}
\end{displaymath}
and
\begin{displaymath}
\xymatrix{H_{n}\big(\Der(\mathcal{M}_{X}, \mathcal{M}_{X};
1)\big)\ar[rr]^-{H((\mathcal{M}_{f})^*)} \ar[d]_-{H(\varepsilon_*)}
& & H_{n}\big(\Der(\mathcal{M}_{Y},
\mathcal{M}_{X}; \mathcal{M}_f)\big)\ar[d]^-{H(\varepsilon_*)} \\
H_{n}\big(\Der(\mathcal{M}_{X}, \Q;
\varepsilon)\big)\ar[rr]^-{H(\widehat{(\mathcal{M}_{f})^*})} & &
H_{n}\big(\Der(\mathcal{M}_{Y}, \Q; \varepsilon)\big)}
\end{displaymath}
are equivalent for each $n\geq2$.
\end{theorem}

We prove this result below.  First, we comment on some ingredients
of the statement and proof, and give some immediate consequences.

\begin{remark}\label{rem:indecomposables}
In rational homotopy theory there is a standard way to identify
the rational homotopy groups of a space $X$, and more generally
the homomorphism of rational homotopy groups induced by a map of
spaces.  Namely, the rational homotopy groups are identified with
the dual of the vector space of indecomposables of the minimal
model, thus $\pi_{*}(X)\otimes\Q \cong \Hom(Q(\mathcal{M}_{X}),
\Q)$.  For maps, $f_\#\otimes\Q$ is identified with the dual of
the map of vector spaces of indecomposables $Q(\mathcal{M}_{f})
\colon Q(\mathcal{M}_{Y}) \to Q(\mathcal{M}_{X})$ induced by the
Sullivan minimal model $\mathcal{M}_{f} \colon \mathcal{M}_{Y} \to
\mathcal{M}_{X}$ of a map $f \colon X \to Y$ (see
\cite[Sec.15(d)]{F-H-T} for details). From the bottom maps in the
two squares of \thmref{thm:equivalence of squares}, we obtain a
superficially different description of $f_\#\otimes\Q$.  But it is
easy to see that this agrees with the standard one: Note that the
derivation law implies $\Der_{*}(\mathcal{M}_{X}, \Q; \varepsilon)
\cong \Hom(Q(\mathcal{M}_{X}), \Q)$,  while the minimality of
$\mathcal{M}_{X}$ implies that $\delta = 0$ in the chain complex
$\Der_{*}(\mathcal{M}_{X}, \Q; \varepsilon)$. Thus we have
$\Hom(Q(\mathcal{M}_{X}), \Q) \cong
H_{*}\big(\Der(\mathcal{M}_{X}, \Q; \varepsilon)\big)$.  The
agreement between our description of $f_\#\otimes\Q$ and the
standard one is obvious from this isomorphism.
\end{remark}

\thmref{thm:equivalence of squares} also contains and depends upon
basic results concerning the rationalization of function space
components and evaluation subgroups due to several authors.  We
consider this material here.

  When $X$ is a finite complex, the function space
 $\map(X, Y)$ has the homotopy type of a CW complex by the result of  Milnor
 \cite{Mil}.
 In fact, by Hilton-Mislin-Roitberg
 \cite[Ch. II, Th. 2.5]{H-M-R} the components $\map(X, Y;f)$ are
 nilpotent complexes.
 Moreover, given a rationalization $e_{Y} \colon Y \rightarrow
Y_{\Q}$ of $Y$ the induced map $(e_{Y})_{*} \colon \map(X,Y;f)
\rightarrow \map(X,Y_{\Q}; e_{Y}\circ f)$ is a rationalization of
$\map(X,Y;f)$ \cite[Ch. II, Th. 3.11]{H-M-R}. By
\cite[Th.2.3]{Sm}, rationalization in the initial variable $e_{X}
\colon X \rightarrow X_{\Q}$ induces a weak equivalence
$(e_{X})^{*} \colon \map(X_{\Q},  Y_{\Q}; f_\Q ) \rightarrow
\map(X, Y_{\Q}; f_\Q\circ e_{X})$.  These results,  together with
the naturality of the various maps involved,  imply the following
result.

\begin{theorem} \label{thm:rationalization} Let $f \colon X \to Y$ be a map between  simply connected
complexes of finite type with $X$ finite. Let $f_{\Q} : X_{\Q} \to
Y_{\Q}$ denote the rationalization of $f$. The commutative squares
\begin{displaymath}
\xymatrix{\pi_n\big(\map(X, X; 1)\big)\otimes\Q
\ar[rr]^{(f_*)_\#\otimes\Q} \ar[d]_{\omega_\#\otimes\Q} & &
\pi_n\big(\map(X, Y;
f)\big)\otimes\Q\ar[d]^{\omega_\#\otimes\Q}\\
\pi_n(X)\otimes\Q \ar[rr]_{f_\#\otimes\Q} & & \pi_n(Y)\otimes\Q}
\end{displaymath}
and
\begin{displaymath}
\xymatrix{\pi_n\big(\map(X_\Q, X_\Q; 1)\big)
\ar[rr]^{((f_\Q)_*)_\#} \ar[d]_{(\omega)_\#} & &
\pi_n\big(\map(X_\Q, Y_\Q;
f_\Q)\big)\ar[d]^{(\omega)_\#}\\
\pi_n(X_\Q) \ar[rr]_{(f_\Q)_\#} & & \pi_n(Y_\Q)}
\end{displaymath}
are equivalent for each $n \geq 2.$
\end{theorem}
The  result  $G_n(X_\Q) \cong G_n(X) \otimes \Q$ for $X$ a  simply
connected  finite complex  due to Lang \cite{Lan}  is an easy
consequence of  \thmref{thm:rationalization}, as is its
generalization $G_*(Y_\Q, X_\Q; f_\Q) \cong G_*(Y, X;f) \otimes
\Q$ (c.f.~\cite{Woo-Kim2, Sm}) for $X, Y$ simply connected
complexes of  finite type with $X$ finite. The corresponding
localization result for the relative Gottlieb group can be deduced
from the preceding discussion, as well.
\begin{theorem} \label{thm:relative Gottlieb rationalization}  Let
    $f \colon X \to Y$ be a map between  simply connected
complexes of finite type with $X$ finite. Then
$G^{rel}_n(Y_\Q, X_\Q; f_\Q) \cong G^{rel}_n(Y, X;f) \otimes \Q,$ for
$n \geq 3.$
\end{theorem}
\begin{proof}
Consider, as in the introduction,  the long exact homotopy sequence of
the induced map  $f_* \colon  \map(X,X; 1) \to \map(X, Y; f)$ and
its relative homotopy group $\pi_n(f_*)$ for $n \geq 3.$   We also
have
$(f_\Q)_* \colon  \map(X_\Q,X_\Q; 1) \to \map(X_\Q, Y_\Q; f_\Q)$
and the   relative group $\pi_n((f_\Q)_*).$      The  results
cited above and the $5$-Lemma, imply  the  maps  induced on the
various function
spaces  involved   by  the rationalizations
$e_X : X \to X_\Q$ and
$e_Y : Y \to Y_\Q$ induce  a rationalization homomorphism
$R \colon   \pi_n(f_*) \to \pi_n((f_\Q)_*).$ The maps $e_X$ and
$e_Y$ also induce  a rationalization $r \colon \pi_n(f) \to \pi_n(f_\Q).$
 By naturality, we obtain a commutative diagram
 \begin{displaymath}
     \xymatrix{   \pi_n(f_*) \ar[r]^{R} \ar[d]^{\omega_\sharp}
     & \pi_n((f_\Q)_*) \ar[d]^{\omega_\sharp} \\
     \pi_n(f) \ar[r]^{r} & \pi_n(f_\Q).
     }
     \end{displaymath} Thus  $
G^{rel}_n(X_\Q, Y_\Q; f_\Q) = \omega_\sharp \left( \pi_n((f_\Q)_*
\right) =
  r \circ \omega_\sharp \left(\pi_n(f_*) \right) \cong
    G_n^{rel}(Y, X ;f) \otimes \Q
    $
     \end{proof}

     Combining the preceding, we obtain the following:

\begin{corollary} \label{cor:rationalization of the G-sequence}
     Let
    $f \colon X \to Y$ be a map between  simply connected
complexes of finite type with $X$ finite.  Then the
rationalization of the $G$-sequence of $f$
$$\cdots \rightarrow G_{n}(X)\otimes\Q \rightarrow G_{n}(Y, X;f)\otimes\Q \rightarrow
G_{n}^{rel}(Y,X;f)\otimes\Q \rightarrow G_{n-1}(X)\otimes\Q
\rightarrow \cdots $$
is equivalent to the $G$-sequence of the rationalization of $f$
$$\cdots \rightarrow G_{n}(X_\Q) \rightarrow G_{n}(Y_\Q, X_\Q;f_\Q) \rightarrow
G_{n}^{rel}(Y_\Q,X_\Q;f_\Q) \rightarrow G_{n-1}(X_\Q) \rightarrow
\cdots
$$
\end{corollary}

In the next section, we identify the rational $G$-sequence  in the
context of rational homotopy theory.  Meanwhile, from
\thmref{thm:equivalence of squares}  we immediately retrieve
minimal model descriptions of  the rationalized Gottlieb groups of
a space, and of the rationalized evaluation subgroups of a map.
The first of these is well-known:

\begin{corollary}\label{cor:rational Gottlieb group}
Let $X$ be a simply connected finite complex.  The rationalized
$n$th Gottlieb group $G_n(X_{\Q})\cong G_n(X)\otimes \Q$ is
isomorphic to the image of the induced homomorphism
$$H(\varepsilon_*)\colon H_{n}\big(\Der(\mathcal{M}_{X}, \mathcal{M}_{X};1)\big) \to H_{n}\big(\Der(\mathcal{M}_{X}, \Q;
\varepsilon)\big)$$
for $n \geq 2.$
\end{corollary}

This is easily translated into the standard minimal model
description of the rationalized Gottlieb groups given by F{\'e}lix
and Halperin (see \cite[Sec.29(c)]{F-H-T}). Specifically, they
describe a Gottlieb element of the minimal model $\mathcal{M}_{X}
= \Lambda(V)$ as a linear map $\theta \colon V^n \to \Q$ that
extends to a derivation of $\mathcal{M}_{X}$ satisfying $d\theta =
(-1)^{n} \theta d$. Such a derivation $\theta$ is a cycle in
$\Der(\mathcal{M}_{X}, \mathcal{M}_{X};1)$, and the class that it
represents has non-zero image under $H(\varepsilon_*)$ precisely
when the original linear map $\theta \colon V^n \to \Q$ is
non-zero. On recalling that $H_{n}\big(\Der(\mathcal{M}_{X}, \Q;
\varepsilon)\big) \cong \Der(\mathcal{M}_{X}, \Q; \varepsilon)
\cong \Hom(Q(\mathcal{M}_{X}), \Q)$, we see the two descriptions
agree.

\begin{corollary}\label{cor:rational eval subgp of f}
Let $f \colon X \to Y$ be a map between  simply connected
complexes of finite type with $X$ finite. The rationalized $n$th
evaluation subgroup $G_n(Y_{\Q}, X_{\Q}; f_{\Q})\cong G_n(Y, X; f)
\otimes {\Q}$ of the map $f$ is isomorphic to the image of the
induced homomorphism
\begin{displaymath}
H(\varepsilon_*)\colon H_{n}\big(\Der(\mathcal{M}_{Y},
\mathcal{M}_{X};\mathcal{M}_{f})\big) \to
H_{n}\big(\Der(\mathcal{M}_{Y}, \Q; \varepsilon)\big).
\end{displaymath}
for $n \geq 2.$
\end{corollary}

This identification of the rationalized evaluation subgroups of
the map $f$ can be conveniently phrased in a way comparable to the
F{\'e}lix-Halperin description of the rationalized Gottlieb
groups: An evaluation subgroup element of the minimal model
$\mathcal{M}_{f} \colon \mathcal{M}_{Y} \to \mathcal{M}_{X}$, with
$\mathcal{M}_{Y} = \Lambda(W)$,  is a linear map $\theta \colon
W^n \to \Q$ that extends to an $\mathcal{M}_{f}$-derivation
$\theta  \colon \mathcal{M}_{Y} \to \mathcal{M}_{X}$ satisfying
$d_Y \theta  = (-1)^{n} \theta d_X$. Whenever such a $\theta \in
\Der_n(\mathcal{M}_{Y}, \Q;\varepsilon)$ is non-zero, it is a
non-zero element in the image of $H(\varepsilon_*)$.

In view of the preceding remarks, we introduce the following
vocabulary and notation.

\begin{definition}\label{def:eval subgp of DG map}
Suppose $\phi \colon A \to B$ is a map of DG algebras.  We define
the \emph{evaluation subgroup of $\phi$} as the image of the map
\begin{displaymath}
H(\varepsilon_*)\colon H_{n}\big(\Der(A, B; \phi)\big) \to
H_{n}\big(\Der(A, \Q; \varepsilon)\big).
\end{displaymath}
We denote it by $G_n(A, B; \phi)$.  In the special case in which
$A = B$ and $\phi = 1_B$, we refer to the \emph{Gottlieb group of
$B$}, and use the notation $G_n(B)$.
\end{definition}

From the previous discussion, we see that $G_n(\mathcal{M}_{Y},
\mathcal{M}_{X}; \mathcal{M}_{f}) \cong G_n(Y, X; f) \otimes \Q$
and $G_n(\mathcal{M}_{X}) \cong G_n(X) \otimes \Q$.

\begin{proof}[Proof of \thmref{thm:equivalence of squares}]
We will define vector space isomorphisms $\Phi$, $\Phi_f$,
$\Psi_X$, and $\Psi_Y$ to give the following equivalence of
commutative squares:
\begin{displaymath}
\xymatrix@C=-32pt{ &
H_{n}\big(\Der(\mathcal{M}_{X},\mathcal{M}_{X};1)\big)
 \ar[rr]^-{H((\mathcal{M}_{f})^*)}
\ar'[d]^-{H(\varepsilon_*)}[dd] & &
H_{n}\big(\Der(\mathcal{M}_{Y},\mathcal{M}_{X};\mathcal{M}_{f})\big)
\ar[dd]^-{H(\varepsilon_*)} \\
\pi_n\big(\map(X, X; 1)\big)\otimes\Q
\ar[ru]^{\Phi}_{\cong}\ar[rr]_(0.5){(f_*)_\#\otimes1}
 \ar[dd]_{\omega_\#}& &
\pi_n\big(\map(X, Y; f)\big)\otimes\Q \ar[ru]^{\Phi_f}_{\cong} \ar[dd]^(0.3){\omega_\#}\\
 &  H_{n}\big(\Der(\mathcal{M}_{X}, \Q;\varepsilon)\big)
\ar'[r]_(0.9){H(\widehat{(\mathcal{M}_{f})^*})}[rr] & &
H_{n}\big(\Der(\mathcal{M}_{Y},
\Q;\varepsilon)\big)  \\
\pi_n(X)\otimes\Q \ar[ru]^{\Psi_X}_{\cong} \ar[rr]_{f_\#\otimes1}
& & \pi_n(Y)\otimes\Q\ar[ru]^{\Psi_Y}_{\cong} }
\end{displaymath}

We obtain $\Phi_{f}$ as the rationalization of a natural
homomorphism
$$\Phi'_{f} \colon \pi_{n}(\map(X,Y;f)) \to
H_{n}\big(\Der(\mathcal{M}_{Y},\mathcal{M}_{X};\mathcal{M}_{f})\big),$$
which we now define.  A representative of a homotopy class $\alpha
\in \pi_{n}(\map(X,Y;f))$  determines, via the exponential
correspondence,  a map $F \colon S^{n} \times X \rightarrow Y$
that satisfies $F \circ i_2 = f$, where $i_2 \colon X \rightarrow
S^{n} \times X$ is the inclusion.  The map $F$ is often called an
\emph{affiliated map} for (a representative of) $\alpha$.  Passing
to minimal models, we obtain a map $\mathcal{M}_F \colon
\mathcal{M}_{Y} \rightarrow \mathcal{M}_{S^{n}} \otimes
\mathcal{M}_{X}$ with $\mathcal{M}_{i_2} \circ \mathcal{M}_F =
\mathcal{M}_f$ (equals, not just up to DG homotopy---see
\propref{prop:restricted DG homotopy}). Now $S^n$ is a
\emph{formal} space, which means there is a quasi-isomorphism of
DG algebras $\psi \colon \mathcal{M}_{S^{n}} \rightarrow
H^{*}(S^{n}; \Q)$.  In turn, this gives a quasi-isomorphism
$\psi\otimes 1 \colon \mathcal{M}_{S^{n}}\otimes \mathcal{M}_{X}
\rightarrow H^{*}(S^{n}; \Q) \otimes \mathcal{M}_{X}$. Recall that
$\Lambda V$ denotes the free algebra generated by $V$.  That is,
polynomial on generators of even degree, and exterior on
generators of odd degree. Write $H^{*}(S^{n}; \Q)$ as
$\Lambda(s_n)/(s_n^2)$ if $n$ is even, or as $\Lambda(s_n)$ if $n$
is odd. Given $\chi \in \mathcal{M}_{Y},$ we may write
$$(\psi\otimes 1)\circ\mathcal{M}_F(\chi) = 1 \otimes \mathcal{M}_f(\chi) + s_{n}
\otimes \theta_G(\chi),$$
thus defining a linear map $\theta_G\colon \mathcal{M}_{Y} \to
\mathcal{M}_{X}$ that reduces degree by $n$.   A standard
check---using the fact that $(\psi\otimes 1)\circ\mathcal{M}_F$ is
a DG algebra map---shows that $\theta_G$ is an
$\mathcal{M}_f$-derivation that is a
$\delta_{\mathcal{M}_f}$-cycle. Set $\Phi'_{f}(\alpha) =
[\theta_G] \in H_{n}\big(\Der(\mathcal{M}_{Y}, \mathcal{M}_{X};
\mathcal{M}_f)\big)$.

To show $\Phi'_{f}(\alpha)$ is well-defined, suppose that $g_1,
g_2 \colon S^n \to \map(X,Y;f)$ are homotopic representatives of
$\alpha$ with affiliated maps $F, G \colon S^{n} \times X
\rightarrow Y$ respectively.   Then the homotopy $K \colon S^n
\times I \to \map(X,Y;f)$ from $g_1$ to $g_2$ gives a homotopy $H
\colon S^n \times X \times I \to Y$, from $F$ to $G$, by setting
$H(s,x,t) = K(s,t)(x)$.   Further, since $K$ is a \emph{based
homotopy}, the homotopy $H$ satisfies $H\circ i = J$, where $i$
denotes the inclusion $i(x,t) = (*,x,t)$ and $J$ denotes the
homotopy $J(x,t) = f(x)$ that is stationary at $f$. A basic result
of rational homotopy theory says that homotopic maps have DG
homotopic Sullivan minimal models (see \cite[Ch.12]{F-H-T}. So the
homotopy $H$ gives a DG homotopy $\mathcal{H} \colon
\mathcal{M}_{Y} \to \mathcal{M}_{S^{n}} \otimes
\mathcal{M}_{X}\otimes\Lambda(t, dt)$ between minimal models for
$F$ and $G$. Translating the restriction on $H$ into minimal model
terms allows us to assume that $\mathcal{H}$ is such that
\begin{equation}\label{eq:DG homotopy K}
(\psi\otimes 1\otimes 1)\circ \mathcal{H}(\chi) =
1\otimes\mathcal{M}_f(\chi)\otimes1 + \sum_{i \geq 0}
s_n\otimes\alpha_i(\chi)\otimes t^i + \sum_{i \geq 0}
s_n\otimes\beta_i(\chi)\otimes t^idt,
\end{equation}
for an element $\chi \in \mathcal{M}_{Y}$. This translation is
intuitively plausible, but its justification requires some
technical details, which we provide in \propref{prop:restricted DG
homotopy}.  Since the DG homotopy $\mathcal{H}$ is from
$\mathcal{M}_F$ to $\mathcal{M}_G$, then at $t = 0$ we have
$\alpha_0(\chi) = \theta_F(\chi)$, and from $t=1$, we have
$\sum_{i \geq 0} \alpha_i(\chi)= \theta_G(\chi)$. To establish
well-definedness, we must show these differ by a boundary in
$\Der(\mathcal{M}_{Y}, \mathcal{M}_{X}; \mathcal{M}_f)$. To this
end, use (\ref{eq:DG homotopy K}) to write separate expressions
for $(\psi\otimes 1\otimes 1)\circ \mathcal{H}(\chi\chi')$ and
$(\psi\otimes 1\otimes 1)\circ \mathcal{H}(\chi)(\psi\otimes
1\otimes 1)\circ \mathcal{H}(\chi')$. Since a DG homotopy is an
algebra map, these expressions agree.  By equating them and
collecting like terms we obtain equations
$$\beta_{i}(\chi \chi') = (-1)^{n|\chi|} \mathcal{M}_f(\chi)\beta_i(\chi') +
(-1)^{|\chi'|} \beta_i(\chi)  \mathcal{M}_f(\chi')
$$
for each $i \geq 0$.  By substituting $\gamma_i(\chi) =
(-1)^{|\chi|}\beta_i(\chi)$ for $\chi \in \mathcal{M}_{Y}$, we
obtain  derivations $\gamma_i \in \Der_{n+1}(\mathcal{M}_{Y},
\mathcal{M}_{X}; \mathcal{M}_f)$. On the other hand, use
(\ref{eq:DG homotopy K}) to write separate expressions for
$(\psi\otimes 1\otimes 1)\circ \mathcal{H}(d\chi)$ and
$d(\psi\otimes 1\otimes 1)\circ \mathcal{H}(\chi)$, with the
latter obtained by applying $d$ to both sides of (\ref{eq:DG
homotopy K}). Since a DG homotopy respects differentials, these
expressions agree. By equating them and collecting like terms we
obtain equations
$$\beta_{i}(d\chi) = (-1)^{|\chi|} (i+1) \alpha_{i+1}(\chi) +
(-1)^{n} d\beta_i(\chi)
$$
for each $i \geq 0$.  With the previous substitution, this gives
$d\gamma_i(\chi) - (-1)^{n+1} \gamma_i d(\chi) = (-1)^{n+1} (i+1)
\alpha_i(\chi)$, that is,
$$\alpha_{i+1}(\chi) = \delta_{\mathcal{M}_f} \big((-1)^{n+1}  \frac{1}{i+1}\,\gamma_i(\chi)\big)$$
for each $i \geq 0$.  It follows that the difference of
derivations $\theta_G - \theta_F = \sum_{i \geq 1} \alpha_i$ is a
$\delta_{\mathcal{M}_f}$-boundary in $\Der(\mathcal{M}_{Y},
\mathcal{M}_{X}; \mathcal{M}_f)$.  Hence $\Phi'_f$ is
well-defined.

It is not difficult to show that $\Phi'_f$ is a homomorphism.  But
since the proof requires some technical notions from rational
homotopy theory, we postpone it to the appendix (\propref{prop:Phi
homomorphism}). As we stated before, the map of vector spaces
$\Phi_f$ is now obtained as the rationalization of the (group)
homomorphism $\Phi'_f$. The map $\Phi$ is defined in the same way,
specializing to the case in which $Y = X$ and $f = 1_X$. The maps
$\Psi_X$ and $\Psi_Y$ are the standard minimal model
identification of the rational homotopy groups of a space, as
discussed in \remref{rem:indecomposables}.

Next, we show $\Phi_{f}$ is surjective.  Denote by
$[S^n_{\Q}\times X_{\Q}, Y_{\Q}]_{f_{\Q}}$ the subset of the set
of homotopy classes of maps $S^n_{\Q}\times X_{\Q} \to Y_{\Q}$
consisting of classes represented by a map that restricts to
$f_{\Q}$ on $X_{\Q}$.  By \thmref{thm:rationalization}, we
identify $\pi_{n}(\map(X,Y;f))\otimes\Q$ with
$\pi_{n}(\map(X_{\Q}, Y_{\Q}; f_{\Q}))$, and hence with
$[S^n_{\Q}\times X_{\Q}, Y_{\Q}]_{f_{\Q}}$. Now suppose given
$[\theta] \in H_{n}\big(\Der(\mathcal{M}_{Y},\mathcal{M}_{X};
\mathcal{M}_f)\big)$. Use $\theta$ to define
\begin{displaymath}
\phi(\chi) = 1 \otimes \mathcal{M}_f(\chi) + s_{n} \otimes
\theta(\chi)
\end{displaymath}
for $\chi \in \mathcal{M}_{Y}$.  Since $\theta$ is an
$\mathcal{M}_f$-derivation that is a cycle, this defines a DG
algebra map $\phi \colon \mathcal{M}_{Y} \to H^{*}(S^{n};\Q)
\otimes \mathcal{M}_{X}$.  Now lift $\phi$ through the surjective
quasi-isomorphism $\psi\otimes1$ as in \cite[Lem.12.4]{F-H-T}, to
obtain a map $\widetilde{\phi} \colon \mathcal{M}_{Y} \rightarrow
\mathcal{M}_{S^{n}} \otimes \mathcal{M}_{X}$ that satisfies
$(\varepsilon\cdot1)\circ\widetilde{\phi} = \mathcal{M}_f \colon
\mathcal{M}_{Y} \to \mathcal{M}_{X}$.  By the standard
correspondence between maps of minimal models and maps of rational
spaces, this gives a map $F \colon S^{n}_{\Q} \times X_{\Q}
\rightarrow Y_{\Q}$ that satisfies $i_2\circ F \sim f_{\Q} \colon
X_{\Q} \rightarrow Y_{\Q}$.  Using, for example,
\cite[Th.9.7]{F-H-T}, we can adjust $F$ into a homotopic map $F'
\colon S^{n}_{\Q} \times X_{\Q} \to Y_{\Q}$ that satisfies
$i_2\circ F' = f_{\Q} \colon X_{\Q} \to Y_{\Q}$, so that $F'$
represents a class of $[S^n_{\Q}\times X_{\Q}, Y_{\Q}]_{f_{\Q}}$.
As described at the start of this paragraph, $F'$ corresponds to a
homotopy class $\alpha \in \pi_{n}(\map(X,Y;f))\otimes\Q$.
Evidently, we have $\Phi_{f}(\alpha) = \theta$.

Finally, we show $\Phi_{f}$ is injective.  Since $\Phi_{f}$ is a
vector space homomorphism, it is sufficient to show that $\alpha
\in  \pi_{n}(\map(X,Y;f))\otimes\Q$ is zero whenever
$\Phi_{f}(\alpha) = 0$.  Using the identification of the previous
paragraph, let $G \colon S^n_{\Q}\times X_{\Q} \to Y_{\Q}$ be an
affiliated map for $\alpha$. Suppose that $\theta_G =
\delta(\eta)$ for $\eta \in
\Der_{n+1}(\mathcal{M}_{Y},\mathcal{M}_{X}; \mathcal{M}_f)$. Using
$\eta$, define a map $\Gamma \colon \mathcal{M}_{Y} \to
H^{*}(S^{n};\Q) \otimes \mathcal{M}_{X}\otimes\Lambda(t, dt)$ by
\begin{displaymath}
\Gamma(\chi) = 1\otimes\mathcal{M}_f(\chi)\otimes1 +
s_n\otimes\theta_G(\chi)\otimes (1-t) +
s_n\otimes\eta(\chi)\otimes dt.
\end{displaymath}
A routine check verifies that $\Gamma$ is a DG algebra map.
Furthermore, it is a DG homotopy from $\mathcal{M}_{G}$ to the map
$E \colon \mathcal{M}_{Y} \to H^{*}(S^{n};\Q) \otimes
\mathcal{M}_{X}$ given by $E(\chi) =
1\otimes\mathcal{M}_{f}(\chi)$.  Now this latter map is a Sullivan
model of an affiliated map for $0 \in
\pi_{n}(\map(X,Y;f))\otimes\Q$. Therefore, the DG homotopy
translates into a homotopy between affiliated maps $S^n_{\Q}\times
X_{\Q} \to Y_{\Q}$ for $\alpha$ and $0$. It follows that
$\alpha=0$, and thus $\Phi_{f}$ is injective.  We observe that,
strictly speaking, we have not justified that the homotopy between
the maps $S^n_{\Q}\times X_{\Q} \to Y_{\Q}$ is relative to
$X_{\Q}$, which corresponds to the homotopy between the maps $S^n
\to \map(X_{\Q}, Y_{\Q};f_{\Q})$ being based. However, a based map
from a sphere is null-homotopic if and only if it is based
null-homotopic (cf. \cite[p.27]{Spa89}).

Commutativity of the cube diagram follows from the naturality of
the homomorphism $\Phi'_{f}$.  By naturality, we mean the
following:  Suppose given maps of spaces $f\colon A \to B$ and
$g\colon B \to C$.  The we have induced maps of function spaces
$g_* \colon \map(A,B;f) \to \map(A,C; g\circ f)$ and $f^* \colon
\map(B,C;g) \to \map(A,C; g\circ f)$.  For either case we obtain a
commutative square involving $\Phi_{g\circ f}$.  Namely, we have
$H\big((\mathcal{M}_g)^*\big) \circ\Phi_{f} = \Phi_{g\circ f}\circ
(g_*)_{\#}$ in the first case, and $H\big((\mathcal{M}_f)_*\big)
\circ\Phi_{g} = \Phi_{g\circ f}\circ (f^*)_{\#}$ in the second
case, as is easily checked.  Since $\Phi_{f}$ is obtained from
$\Phi'_{f}$ by localization, the isomorphism $\Phi_{f}$ has the
same naturality property.  This is sufficient to conclude that the
top, bottom, left, and right faces of the cube commute.  For the
evaluation map $\omega \colon \map(X,Y;f) \to Y$ can be identified
with $i^* \colon \map(X,Y;f) \to \map(x_0,Y;y_0)$, where $x_0\in
X$ and $y_0 \in Y$ denote basepoints, and $i\colon x_0 \to X$
inclusion of the basepoint. Likewise for the evaluation map
$\omega \colon \map(X,X;1) \to X$, and then $f \colon X \to Y$ can
be identified with $f_* \colon \map(x_0,X;x_0) \to
\map(x_0,Y;y_0)$ (cf.~also \remref{rem:indecomposables}). Finally,
the front and rear faces commute because the squares
(\ref{eq:homotopy square}) and (\ref{eq:derivation square}) are
commutative.
\end{proof}

We finish this section with two more immediate consequences of
\thmref{thm:equivalence of squares}.  The first retrieves a basic
result of Thom, in the rational homotopy setting.

\begin{corollary}(\cite[Th.2]{Th})\label{cor:Thom}
Let $Y = K(V, m)$ be an Eilenberg-Mac Lane space, for $V$ a finite
dimensional (ungraded) rational vector space.   If $X$ is a finite
CW complex, and $f\colon X \to Y$ is any map, then
$\pi_n\big(\map(X,Y;f)\big) \cong H^{m-n}(X;V)$.
\end{corollary}

\begin{proof}
The minimal model for $Y$ is $\Lambda V^*$ with zero differential,
where $V^*$ denotes  the dual vector space of $V$. It follows
easily that
$H_{n}\big(\Der(\mathcal{M}_{Y},\mathcal{M}_{X};\mathcal{M}_{f})\big)
\cong \Hom(V^*, H^{m-n}(\mathcal{M}_{X})) \cong \Hom(V^*,
H^{m-n}(X;\Q)) \cong H^{m-n}(X;V)$.
\end{proof}

Notice that this result---with the remarks on rationalization
preceding \thmref{thm:rationalization}---easily extends to yield
the rational homotopy type of $\map(X,Y;f)$, in case $Y$ is a
rational $H$-space. That is, a space whose rationalization is an
$H$-space.  For if $Y$ is a rational $H$-space, then so too is
$\map(X,Y;0)$, which is homotopy equivalent---via translation by
$f$---to $\map(X,Y;f)$. Now a rational $H$-space is determined up
to rational homotopy type by its rational homotopy groups.
Furthermore, a rational $H$-space is a product of rational
Eilenberg-Mac Lane spaces. Hence $\map(X,Y;f)$ has the rational
homotopy type of a product of spaces $\map(X,Y_i;0)$, for $Y_i$ a
rational Eilenberg-Mac Lane space as in the corollary.

We give a  further consequence of \thmref{thm:equivalence of
squares} concerning the rational homotopy groups of certain
function spaces. Define an \emph{$F_{0}$-space} to be a finite
simply connected complex with finite dimensional rational homotopy
(a \emph{rationally elliptic space}) such that $H^{odd}(X,\Q) =
0$. This type of space features in the following well-known
conjecture of Halperin (cf.~\cite[p.516]{F-H-T}):

\begin{conjecture}\label{conj:Halperin}
Suppose $X$ is an $F_0$-space.  Then any fibration $X \to E \to B$
of simply connected spaces is TNCZ, that is, the fibre inclusion
$j \colon X \to E$ induces a surjection on rational cohomology.
\end{conjecture}

The rational homotopy of $\map(X,X; 1)$ for $X$ an $F_{0}$-space
is directly related to \conjref{conj:Halperin} (see \cite{Me} for
details). The following result extends \cite[Prop.2.6]{Me} (also
compare \cite[Cor.4.6]{Griv} and \cite[Th.B]{Hau}):

\begin{corollary}\label{cor:Grivel}
Let $f\colon X \rightarrow Y$ be a map between $F_{0}$-spaces.
Then for $r \geq 1$,
$$\pi_{2r}(\map(X,Y; f)) \otimes \Q \cong \Der_{2r}(H^{*}(Y, \Q),
H^{*}(X,\Q); H(f)).$$
\end{corollary}

\begin{proof} The argument given by Grivel---for the case in which $Y = X$
and $f = 1_X$---can be used word for word to show
$$H_{2r}(\Der(\mathcal{M}_{Y}, \mathcal{M}_{X}; \mathcal{M}_f))
\cong \Der_{2r}(H^{*}(Y, \Q), H^{*}(X, \Q); H(f)).$$
The result now follows from \thmref{thm:equivalence of squares}.
\end{proof}

\section{Derivation Spaces and Long Exact Sequences}%
\label{sec:long exact sequences}

In this section, we identify the rationalized long exact homotopy
sequences of the maps $f \colon X \to Y$ and $f_* \colon \map(X,X;
1) \to \map(X,Y; f)$, and hence the rationalized $G$-sequence of
$f$.  Our identifications flow from the following observation:
Suppose given vector space homomorphisms $\phi_n \colon A_n \to
B_n$ for each $n$.  Then, up to equivalence, \emph{there is a
unique way to fit these into a three-term long exact sequence}.
Since we rely on this observation for our basic results, we make a
formal statement of the fact.

\begin{lemma}\label{lem:unique les}
Suppose given equivalences of vector space homomorphisms
\begin{displaymath}
\xymatrix{A_n \ar[r]^{\phi_n} \ar[d]^{\alpha_n}_{\cong} & B_n
\ar[d]^{\beta_n}_{\cong}\\
A'_n \ar[r]_{\phi'_n} & B'_n}
\end{displaymath}
for each $n$.  Then for any two long exact sequences of vector
spaces containing the $\phi_n$ and the $\phi'_n$ thus
$$
\xymatrix{\cdots \ar[r] &
C_{n+1}\ar@{.>}[d]^{\gamma_{n+1}}_{\cong} \ar[r]^{\partial_{n+1}}&
A_n \ar[r]^{\phi_n}\ar[d]^{\alpha_n}_{\cong}
& B_n \ar[r]^{p_n}\ar[d]^{\beta_n}_{\cong} & C_n \ar[r]\ar@{.>}[d]^{\gamma_n}_{\cong}& \cdots\\
\cdots \ar[r] & C'_{n+1} \ar[r]^{\partial'_{n+1}} & A'_n
\ar[r]_{\phi'_n} & B'_n \ar[r]^{p'_n}& C'_n \ar[r]& \cdots,}
$$
there exist isomorphisms $\gamma_n \colon C_n \to C'_n$ that make
the long exact sequences equivalent.
\end{lemma}

\begin{proof}
We define a map $\gamma_n$ as follows:  Decompose $C_n$ and $C'_n$
as $C_n \cong \mathrm{im}(p_n)\oplus D_n$ and $C'_n \cong
\mathrm{im}(p'_n)\oplus D'_n$, where $D_n$ and $D'_n$ are
complements.  For $x = p_n(y) \in \mathrm{im}(p_n)$, define
$\gamma_n(x) = g'_n\circ\beta_n(y)$.  For $x \in D_n$, set
$\gamma_n(x) = y'$, where $y' \in D'_n$ is such that
$\partial'_n(y'_n) = \alpha_{n-1}\circ\partial_n(x)$.  It is
straightforward to check that $\gamma_n$ is a well-defined
isomorphism.  Indeed, it restricts to give isomorphisms
$\mathrm{im}(p_n)\cong \mathrm{im}(p'_n)$ and $D_n \cong D'_n$.
The required commutativity properties follow immediately.
\end{proof}

There is nothing remarkable in this observation.  It is important,
nonetheless, since it allows us to choose descriptions of the long
exact homotopy sequences that we need in whatever way is most
convenient for our purposes.

We begin with the long exact homotopy sequence of the map $f_*
\colon \map(X,X; 1)$ $\to \map(X,Y; f)$.  From
\thmref{thm:equivalence of squares}, we see that the map this
induces on rational homotopy groups can be identified with $H\big(
(\mathcal{M}_f)^{*}\big)$. Since this map is a homomorphism
induced on homology by a map of chain complexes, there is a
standard way to fit it into a long exact sequence, which we now
describe.

Given a map $\phi \colon A \to B$ of chain complexes, define a
\emph{relative} chain complex $\Rel_{*}(\phi)$ as follows:
$\Rel_{n}(\phi) = A_{n-1}\oplus B_{n}$, with differential
$\delta_{\phi}$ of degree $-1$ given by $\delta_{\phi}(a , b) =
\big(\delta_{A}(a), \delta_{B}(b) - \phi(a)\big)$. Now define
chain maps $J \colon B_{n}\to \Rel_{n}(\phi)$ and $P \colon
\Rel_{n}(\phi) \to A_{n-1}$ by $J(b) = (0, b)$ and $P(a,b) = a$.
On passing to homology, we obtain a long exact sequence of the
following form:
\begin{displaymath}
\xymatrix{\cdots \ar[r] &  H_{n+1}(\Rel(\phi)) \ar[r]^-{H(P)} &
H_{n}(A) \ar[r]^-{H(\phi)} & H_{n}(B) \ar[r]^-{H(J)}&
H_{n}(\Rel(\phi)) \ar[r] & \cdots,}
\end{displaymath}
in which $H_{n}(\Rel(\phi))$ denotes the $n$th homology of the
chain complex $\Rel_{*}(\phi)$. We refer to this long exact
sequence as the \emph{long exact homology sequence of the map
$\phi$}.  Furthermore, this construction is natural. For suppose
given a commutative square of DG vector spaces
\begin{equation}
\xymatrix{A \ar[r]^{\phi}\ar[d]_{\alpha} & B \ar[d]_{\beta}\\
A' \ar[r]^{\phi'} & B'.}
\end{equation}
Then the obvious map $(\alpha, \beta) \colon \Rel_{*}(\phi) \to
\Rel_{*}(\phi')$ is a chain map that satisfies $(\alpha, \beta)J =
J'\beta$ and $\alpha P = P' (\alpha, \beta)$.  Thus we obtain a
homology ladder
\begin{displaymath}
\xymatrix{\cdots \ar[r] &  H_{n+1}(\Rel(\phi)) \ar[r]^-{H(P)}
\ar[d]^-{H(\alpha, \beta)}& H_{n}(A)
\ar[r]^-{H(\phi)}\ar[d]^-{H(\alpha)} & H_{n}(B)
\ar[r]^-{H(J)}\ar[d]^-{H(\beta)}& H_{n}(\Rel(\phi)) \ar[r]\ar[d]^-{H(\alpha, \beta)} & \cdots\\
\cdots \ar[r] &  H_{n+1}(\Rel(\phi')) \ar[r]^-{H(P')} & H_{n}(A')
\ar[r]^-{H(\phi')} & H_{n}(B') \ar[r]^-{H(J')}& H_{n}(\Rel(\phi'))
\ar[r] & \cdots}
\end{displaymath}

In particular, we can apply this construction to the map of chain
complexes
\begin{equation}\label{eq:Mf^* 1}
(\mathcal{M}_{f})^* \colon \Der(\mathcal{M}_{X}, \mathcal{M}_{X};
1) \rightarrow \Der(\mathcal{M}_{Y}, \mathcal{M}_{X};
\mathcal{M}_{f})
\end{equation}
induced by the minimal model $\mathcal{M}_{f} \colon
\mathcal{M}_{Y} \to \mathcal{M}_{X}$ of the map $f \colon X \to
Y$.

\begin{theorem}\label{thm:long exact seq of f_*}
The long exact sequence induced by
$$f_* \colon \map(X,X; 1) \to \map(X,Y; f)$$
on rational homotopy groups is equivalent to the long exact
homology sequence of the map (\ref{eq:Mf^* 1}). Specifically, this
is a long exact sequence
\begin{displaymath}
\xymatrix@C=13pt{& & \cdots \ar[r]^-{H(J)} &
H_{n+1}(\Rel((\mathcal{M}_{f})^*))
 \ar `d[l]  `[llld]_(0.7){H(P)} [llld]
\\
H_{n}\big(\Der(\mathcal{M}_{X}, \mathcal{M}_{X}; 1)\big)
\ar[rr]^-{H((\mathcal{M}_{f})^*)} & &
H_{n}\big(\Der(\mathcal{M}_{Y}, \mathcal{M}_{X};
\mathcal{M}_{f})\big) \ar[r]^-{H(J)} &
\cdots \\
& &  \cdots \ar[r]^-{H(J)} & H_{3}(\Rel((\mathcal{M}_{f})^*))
 \ar `d[l]  `[llld]_(0.7){H(P)} [llld]
\\
 H_{2}\big(\Der(\mathcal{M}_{X},
\mathcal{M}_{X}; 1)\big) \ar[rr]^-{H((\mathcal{M}_{f})^*)} & &
H_{2}\big(\Der(\mathcal{M}_{Y}, \mathcal{M}_{X};
\mathcal{M}_{f})\big) }
\end{displaymath}
in which $\Rel_*((\mathcal{M}_{f})^*)$ is the relative chain
complex of the map (\ref{eq:Mf^* 1}), as described above.
\end{theorem}

\begin{proof}
By \thmref{thm:rationalization} we may assume $f \colon X \to Y$ is
map of rational spaces.   \thmref{thm:equivalence of squares} gives equivalences of vector
space maps
$$\xymatrix{\pi_{n}\big(\map(X,X;1)\big) \ar[d]_{\Phi_{1}}^{\cong}
\ar[rr]^-{(f_*)_\#} & &
\pi_{n}\big(\map(X,Y;f)\big) \ar[d]^{\Phi_{f}}_{\cong}\\
H_{n}\big(\Der(\mathcal{M}_{X}, \mathcal{M}_{X}; 1)\big)
\ar[rr]_-{H\big( (\mathcal{M}_f)^{*}\big)} & &
H_{n}\big(\Der(\mathcal{M}_{Y}, \mathcal{M}_{X};
\mathcal{M}_f)\big),}$$
for each $n \geq 2$.    The top horizontal maps are contained in
the long exact sequence induced by $f_* \colon \map(X,X; 1) \to
\map(X,Y; f)$ on rational homotopy groups.  The bottom horizontal
maps are contained in the long exact homology sequence of the map
(\ref{eq:Mf^* 1}).  From \lemref{lem:unique les}, these sequences
are equivalent.
\end{proof}

\begin{remark}
When we refer to the long exact homotopy sequence of a map, we
mean this in the sense of \cite[Chaps.3,4]{Hilton}:  Recall that
given a map $f \colon X \to Y$, this sequence is as follows:
\begin{displaymath}
     \cdots \rightarrow \pi_{n}(X) \stackrel{f_{\#}}{\longrightarrow}
     \pi_{n}(Y) \rightarrow
\pi_{n}(f) \rightarrow \pi_{n-1}(X) \rightarrow \cdots \rightarrow
\pi_{2}(X) \stackrel{f_{\#}}{\longrightarrow} \pi_{2}(Y).
\end{displaymath}
If $f$ is the inclusion of a subspace, then the groups
$\pi_{n}(f)$ are just the usual homotopy groups of a pair.
Generally, $\pi_{n}(f)$ is defined as homotopy classes of pairs
$(g_{1}, g_{2})$ such that the diagram
$$ \xymatrix{S^{n-1} \ar@{^{(}->}[d] \ar[r]^-{g_{1}} & X \ar[d]^{f} \\
  CS^{n-1} \ar[r]^-{g_{2}}  & Y}
  $$
commutes.  Since $\pi_{2}(f)$ is not necessarily abelian and we
are interested in rationalizing this sequence, we stop at
$\pi_{2}(Y)$.  On the other hand, one can convert $f$ into a
fibration and use the corresponding long exact sequence in
homotopy.  Either approach suits our purposes and indeed the same
sequence of homotopy groups results from either.  From the above,
we see that if $\mathcal{F}$ denotes the homotopy fibre of the map
$f_* \colon \map(X,X; 1) \to \map(X,Y; f)$, then for $n \geq 2$ we
have
\begin{displaymath}
\pi_{n+1}(f^*)\otimes\Q \cong \pi_{n}(\mathcal{F})\otimes\Q \cong
H_{n+1}(\Rel((\mathcal{M}_{f})^*)),
\end{displaymath}
where $\Rel_*((\mathcal{M}_{f})^*)$ is the relative chain complex
of the map (\ref{eq:Mf^* 1}).
\end{remark}

The preceding result specializes to give a description of the long
exact sequence induced by a general map on rational homotopy
groups. The minimal model $\mathcal{M}_{f} \colon \mathcal{M}_{Y}
\to \mathcal{M}_{X}$ of the map $f \colon X \to Y$ also induces a
map of chain complexes
\begin{equation}\label{eq:Mf^* 2}
(\mathcal{M}_{f})^* \colon \Der(\mathcal{M}_{X}, \Q; \varepsilon)
\rightarrow \Der(\mathcal{M}_{Y}, \Q; \varepsilon).
\end{equation}

\begin{theorem}\label{thm:long exact seq of f}
The long exact sequence induced by $f \colon X \to Y$ on rational
homotopy groups is equivalent to the long exact homology sequence
of the map (\ref{eq:Mf^* 2}). Specifically, this is a long exact
sequence
\begin{displaymath}
\xymatrix{  & & \cdots \ar[r]^-{H(J)} &
H_{n+1}(\Rel((\mathcal{M}_{f})^*))
 \ar `d[l]  `[llld]_(0.7){H(P)} [llld]
\\
H_{n}\big(\Der(\mathcal{M}_{X}, \Q;\varepsilon)\big)
\ar[rr]^-{H((\mathcal{M}_{f})^*)} & &
H_{n}\big(\Der(\mathcal{M}_{Y}, \Q;\varepsilon)\big)
\ar[r]^-{H(J)} &
\cdots \\
& &  \cdots \ar[r]^-{H(J)} & H_{3}(\Rel((\mathcal{M}_{f})^*))
 \ar `d[l]  `[llld]_(0.7){H(P)} [llld]
\\
 H_{2}\big(\Der(\mathcal{M}_{X}, \Q;\varepsilon)\big) \ar[rr]^-{H((\mathcal{M}_{f})^*)}
 & &
H_{2}\big(\Der(\mathcal{M}_{Y}, \Q;\varepsilon)\big) }
\end{displaymath}
in which $\Rel_*((\mathcal{M}_{f})^*)$ is the relative chain
complex of the map (\ref{eq:Mf^* 2}).
\end{theorem}

\begin{proof}
Argue exactly as in the proof of \thmref{thm:long exact seq of
f_*}, making use of \thmref{thm:equivalence of squares} and
\lemref{lem:unique les}.
\end{proof}

\begin{remark}
There is already a standard way to describe the long exact
sequence induced by a map on rational homotopy groups, using
minimal models.  This uses the notion of a so-called K-S model of
the map $\mathcal{M}_{f} \colon \mathcal{M}_{Y} \to
\mathcal{M}_{X}$ \cite[Sec.15(d)]{F-H-T}. The description we give
above, however, is better suited to our purposes. Note that if $F$
denotes the homotopy fibre of the map $f \colon X \to Y$, then for
$n \geq 2$ we have
\begin{displaymath}
\pi_{n+1}(f)\otimes\Q \cong \pi_{n}(F)\otimes\Q \cong
H_{n+1}(\Rel((\mathcal{M}_{f})^*)),
\end{displaymath}
where $\Rel_*((\mathcal{M}_{f})^*)$ is the relative chain complex
of the map (\ref{eq:Mf^* 2}).  It is perhaps interesting to
compare the description given in \thmref{thm:long exact seq of f}
to the standard description of the long exact sequence in rational
homotopy groups of a fibration.
\end{remark}

We now turn our attention to the $G$-sequence, and identify it
within our current framework.

Suppose given a DG algebra map $\phi \colon A \to B$.  Starting
from this map, we can construct the following commutative square
of DG vector spaces:
\begin{displaymath}
\xymatrix{\Der_*(B, B; 1)\ar[r]^-{\phi^*} \ar[d]_-{\varepsilon_*}
& \Der_*(A,
B; \phi)\ar[d]^-{\varepsilon_*} \\
\Der_*(B, \Q; \varepsilon)\ar[r]^-{\widehat{\phi^*}} & \Der_*(A,
\Q; \varepsilon).}
\end{displaymath}
In this diagram, $\varepsilon$ denotes the augmentation of either
$A$ or $B$, and we have used a decoration to distinguish the lower
horizontal map from the upper.  On passing to homology and using
the naturality of the relative chain complex construction, we
obtain the following homology ladder:
\begin{displaymath}
\xymatrix{\cdots \ar[r]^-{H(J)} &  H_{n+1}(\Rel(\phi^*))
\ar[r]^-{H(P)} \ar[d]^-{H(\varepsilon_*, \varepsilon_*)}&
H_{n}\big(\Der(B, B; 1)\big)
\ar[r]^-{H(\phi^*)}\ar[d]^-{H(\varepsilon_*)} & H_{n}\big(\Der(A,
B; \phi)\big)
\ar[d]^-{H(\varepsilon_*)} \ \cdots\\
\cdots \ar[r]^-{H(J)} &  H_{n+1}(\Rel(\widehat{\phi^*}))
\ar[r]^-{H(\widehat{P})} & H_{n}\big(\Der(B, \Q; \varepsilon)\big)
\ar[r]^-{H(\widehat{\phi^*})} & H_{n}\big(\Der(A, \Q;
\varepsilon)\big)\  \cdots}
\end{displaymath}
for $n \geq 2$. We supplement \defref{def:eval subgp of DG map}
with the following vocabulary.

\begin{definition}\label{def:rel eval subgp of DG map}
Suppose $\phi \colon A \to B$ is a map of DG algebras.  For $n
\geq 3$ we define the \emph{$n$th relative evaluation subgroup of
$\phi$} as the image of the map
\begin{displaymath}
H(\varepsilon_*, \varepsilon_*)\colon H_{n}\big(\Rel(\phi^*)\big)
\rightarrow H_{n}\big(\Rel(\widehat{\phi^*})\big).
\end{displaymath}
We denote it by $G^{rel}_{n}(A, B; \phi)$.  Then the image of the
upper long exact sequence in the lower, of the ladder above, gives
a (not necessarily exact) sequence
\begin{displaymath}
\xymatrix{ & \cdots \ar[r]^-{H(\widehat{J})} & G^{rel}_{n+1}(A, B;
\phi)
 \ar `d[l]  `[lld]_(0.7){H(\widehat{P})} [lld]
\\
G_n(B) \ar[r]^-{H(\widehat{\phi^*})} & G_n(A, B; \phi)
\ar[r]^-{H(\widehat{J})} & G^{rel}_{n}(A, B; \phi)  \ar`d[l]
`[lld]_(0.7){H(\widehat{P})} [lld]
 \\
\cdots &  \cdots \ar[r]^-{H(\widehat{J})} & G^{rel}_{3}(A, B;
\phi)
 \ar `d[l]  `[lld]_(0.7){H(\widehat{P})} [lld]
\\
 G_{2}(B) \ar[r]^-{H(\widehat{\phi^*})} &
G_{2}(A, B; \phi)\big) }
\end{displaymath}
We refer to this sequence as \emph{the $G$-sequence of the map}
$\phi \colon A \to B$.
\end{definition}

All of the above can be applied to the minimal model
$\mathcal{M}_{f} \colon \mathcal{M}_{Y} \to \mathcal{M}_{X}$ of
the map $f \colon X \to Y$.  By doing so, and then collecting
together previous results, we obtain the following result.

\begin{theorem}\label{thm:G-seq of f}
The rationalization of the $G$-sequence of the map $f \colon X \to
Y$, as far as the term $G_{2}(Y,X;f)$,  is equivalent to the
$G$-sequence of the corresponding map of Sullivan models.
\end{theorem}

\begin{proof}
Starting from the cube displayed in the proof of
\thmref{thm:equivalence of squares}, we extend each of the four
left-to-right maps into their respective long exact sequences.
This is then completed into an equivalence of ladders, by defining
isomorphisms $\gamma_n$ and $\widehat{\gamma_n}$ to give a
commutative square
\begin{displaymath}
\xymatrix{\pi_n(f_*) \ar[r]^-{\gamma_n} \ar[d]_-{\omega_\#} &
H_{n}\big(\Rel(\phi^*)\big)\ar[d]^-{H(\varepsilon_*, \varepsilon_*)} \\
\pi_n(f)\ar[r]^-{\widehat{\gamma_n}} &
H_{n}\big(\Rel(\widehat{\phi^*})\big)}
\end{displaymath}
for each $n \geq 3$.  These isomorphisms are defined as in
\lemref{lem:unique les}, using the top and bottom faces of the
cube.  The one commutativity relation that needs checking, namely
that of the displayed square, follows easily from the
commutativity of the adjacent and parallel squares, together with
the way in which the $\gamma_n$ and $\widehat{\gamma_n}$ are
defined.

The result now follows, since whenever one has such an equivalence
of ladders, the equivalence restricts to give an equivalence of
the corresponding sequences of images.
\end{proof}

In particular, we obtain the companion result to \corref{cor:rational
Gottlieb group} and \corref{cor:rational eval subgp of f}.

\begin{corollary}\label{cor:relative eval subgp of f}
Let $f \colon X \to Y$ be a map between  simply connected
complexes of finite type with $X$ finite.  The rationalized $n$th
relative evaluation subgroup $G^{rel}_n(Y_{\Q}, X_{\Q};
f_{\Q})\cong G^{rel}_n(Y, X; f) \otimes {\Q}$ of the map $f$ is
isomorphic to the image of the induced homomorphism
\begin{displaymath}
H(\varepsilon_*, \varepsilon_*)\colon H_{n}\big(\Rel((\mathcal{M}_f)^*)\big) \to
H_{n}\big(\Rel(\widehat{\mathcal{M}_f^*})\big)
\end{displaymath}
for $n \geq 3.$
\end{corollary}

\begin{remark}
We comment on the low-end terms in the $G$-sequence.  In
\thmref{thm:long exact seq of f} and \thmref{thm:long exact seq of
f_*} we terminate our long exact sequences at the terms
corresponding to $\pi_{2}\big(\map(X, Y;f)\big)$ and $\pi_{2}(Y)$
respectively.  This is because we need a simply connected
hypothesis to ensure our combination of rationalization and
minimal model techniques remains valid.  As a result, our
algebraic description of the rationalized $G$-sequence terminates
at the term corresponding to $G_{2}(X, Y;f)$.    Now in
\thmref{thm:long exact seq of f}, we require $X$ to be simply
connected and finite. As is well-known, this implies
$G_{2i}(X)\otimes\Q = 0$ for each $i$. Therefore, under our
hypotheses, the rationalized $G$-sequence of a map $f \colon X \to
Y$ should be considered as $5$-term (not necessarily exact)
sequences
\begin{displaymath}
\xymatrix{0 \ar[r] & G_{2n}(Y, X;f)\otimes\Q \ar[r] &
G_{2n}^{rel}(Y,X;f)\otimes\Q
 \ar `d[l]  `[lld] [lld]
\\
G_{2n-1}(X)\otimes\Q \ar[r] & G_{2n-1}(Y, X;f)\otimes\Q \ar[r] &
 G_{2n-1}^{rel}(Y,X;f)\otimes\Q \ar[r] & 0
 }
\end{displaymath}
for $n \geq 2$.   Our algebraic description of the rationalized
$G$-sequence given by \thmref{thm:G-seq of f} includes all these
$5$-term sequences.  The ``sporadic" low-end term $G_{2}(Y,
X;f)\otimes\Q$ is best computed by using its characterization
given in \corref{cor:rational eval subgp of f}.
\end{remark}

Before turning to some applications of our algebraic description
of the rationalized $G$-sequence, we give one more description of
a long exact homotopy sequence. Whilst not strictly necessary for
our purposes, it is nonetheless interesting.

We will use $\map_{*}(X,Y; f)$ to denote the \emph{based} mapping
space component. Then we have the fibration sequence
\begin{equation}\label{eq:eval fibn}
\xymatrix{\map_{*}(X,Y;f)\ar[r] &  \map(X,Y; f) \ar[r]^-{\omega} &
Y.}
\end{equation}
We will describe the long exact sequence on rational homotopy
groups induced by this fibration.

Recall that we have the augmentation $\varepsilon \colon A \to \Q$
for a DG algebra $A$.  Let $\widetilde{A}$ denote the
\emph{augmentation ideal}, that is, the kernel of $\varepsilon$.
Given a DG algebra map $\phi \colon A \rightarrow B$, let
$\widetilde{\phi} \colon A \rightarrow \widetilde{B}$ be the DG
algebra map which agrees with $\phi$ in positive degrees and
vanishes in degree zero.

A DG algebra map $\phi \colon A \rightarrow B$ together with the
short exact augmentation sequence
$$\xymatrix{0 \ar[r] & \widetilde{B} \ar[r]^{i} & B
\ar[r]^{\varepsilon} & \Q \ar[r] & 0}$$
of DG algebras gives rise to the short exact sequence of DG vector
spaces
$$
\xymatrix{0 \ar[r] & \Der_{*}(A, \widetilde{B}; \widetilde{\phi})
\ar[r]^{i_*} & \Der_{*}(A, B; \phi) \ar[r]^{\varepsilon_*} &
\Der_{*}(A, \Q; \varepsilon) \ar[r] & 0.}$$
This in turn gives a long exact sequence on homology, in the usual
way, of the form
\begin{equation}\label{eq:augment1}
\xymatrix{ \cdots \ar[r]^-{\Delta_*}& H_{n}\big(\Der(A,
\widetilde{B};\widetilde{\phi})\big) \ar[r]^-{H(i_*)} &
H_{n}\big(\Der(A, B; \phi)\big) \ar[r]^-{H(\varepsilon_*)} &
H_{n}\big(\Der(A, \Q; \varepsilon)\big)\cdots }
\end{equation}
for $n \geq 2$. Call this sequence the \emph{long exact derivation
homology sequence} of the DG algebra map $\phi \colon A
\rightarrow B$.

\begin{theorem}\label{thm:long exact seq of eval fib}
The long exact sequence induced by the evaluation fibration
(\ref{eq:eval fibn}) on rational homotopy groups is equivalent to
the long exact derivation homology sequence (\ref{eq:augment1}) of
the map $\mathcal{M}_{f}\colon \mathcal{M}_{Y} \to
\mathcal{M}_{X}$. Specifically, this is a long exact sequence
\begin{displaymath}
\xymatrix@C=11pt{ & \cdots \ar[r]^-{H(\varepsilon_*)} &
H_{n+1}\big(\Der(\mathcal{M}_{Y}, \Q; \varepsilon)\big)
 \ar `d[l]  `[lld]_(0.7){\Delta_*} [lld]
\\
H_{n}\big(\Der(\mathcal{M}_{Y},
\widetilde{\mathcal{M}_{X}};\widetilde{\mathcal{M}_{f}})\big)
\ar[r]^-{H(i_*)} & H_{n}\big(\Der(\mathcal{M}_{Y},
\mathcal{M}_{X}; \mathcal{M}_{f})\big) \ar[r]^-{H(\varepsilon_*)}
&
\cdots \\
 &  \cdots \ar[r]^-{H(\varepsilon_*)} & H_{3}\big(\Der(\mathcal{M}_{Y}, \Q;
\varepsilon)\big)
 \ar `d[l]  `[lld]_(0.7){\Delta_*} [lld]
\\
H_{2}\big(\Der(\mathcal{M}_{Y},
\widetilde{\mathcal{M}_{X}};\widetilde{\mathcal{M}_{f}})\big)
\ar[r]^-{H(i_*)} & H_{2}\big(\Der(\mathcal{M}_{Y},
\mathcal{M}_{X}; \mathcal{M}_{f})\big) \ar[r]^-{H(\varepsilon_*)}
& H_{2}\big(\Der(\mathcal{M}_{Y}, \Q; \varepsilon)\big).}
\end{displaymath}
\end{theorem}

\begin{corollary}
Let $X$ and $Y$ be simply connected spaces with $X$ finite.  Then
the rational homotopy groups of the based function space
$\map_{*}(X,Y;f)$ are given by
\begin{displaymath}
\pi_{n}\big(\map_{*}(X,Y;f)\big) \otimes\Q \cong
H_{n}\big(\Der(\mathcal{M}_{Y},
\widetilde{\mathcal{M}_{X}};\widetilde{\mathcal{M}_{f}})\big)
\end{displaymath}
for $n \geq 2$.
\end{corollary}

\section{Examples, Computations, and Further Consequences}%
\label{sec:Examples}

We illustrate the effectiveness of the framework established in
the previous two sections with examples. First, we give a
composite example that includes specific computation of many of
the ingredients of the above.  Our example is one in which the
$G$-sequence of a map fails to be exact (after rationalization) at
each of the three types of term that occur.

We begin with some notational conventions. Suppose that
$(A,d_{A})$ and $(B,d_{B})$ are minimal algebras, with $A =
\Lambda(W)$ and $B = \Lambda(V)$ for suitable graded vector spaces
$W$ and $V$. Let  $\phi \colon A \rightarrow B$ be a fixed DG
algebra map.   Since any linear map $W \to B$ extends in a unique
way to a $\phi$-derivation, we can view the space $\Hom_{*}(W, B)$
of negative degree linear maps as a subspace of $\Der_{*}(A,
B;\phi)$. (Of course, in degree $1$ we must restrict to cycles!)
With this view $\Hom_{*}(W, B) = \Der_{*}(A, B;\phi)$  as graded
spaces although the differential depends on the derivation
structure.  This point of view, whereby a derivation is specified
on generators and then extended to the whole algebra, is one that
we will invariably adopt in any practical calculation.

Now suppose given a basis $\{ w_{1},w_{2},w_{3}, \ldots \}$ for
$W$ and an element $P \in B$ with $|P| < |w_{i}|$ we will write $P
\partial w_{i}$ for the $\phi$-derivation carrying $w_{i}$ to $P$
and vanishing on the other $w_{j}$.  Thus any derivation can be
expressed as a sum $\sum_i P_i\partial w_{i}$.   When $B = \Q,$ we
write $w_{j}^{*}$ rather than $1\partial w_{j}$ for the derivation
dual to $w_{j}$.

\begin{example}\label{ex:non-exact G-seq}
Let $f = (f_1, f_2)\colon \H P^{2} \rightarrow S^8\times\H P^{4}$
be the map with coordinate functions $f_1\colon \H P^{2} \to S^8$
obtained by pinching out the bottom cell and $f_2\colon \H P^{2}
\to \H P^{4}$ the inclusion. We will use our framework from above
to compute various terms from the long exact sequences
corresponding to \thmref{thm:long exact seq of f_*} and
\thmref{thm:long exact seq of f}.  Denote $\H P^{2}$ by $X$ and
$S^{8}\times\H P^{4}$ by $Y$, thus $f \colon X \to Y$. Our
computation will show, using \thmref{thm:G-seq of f}, that the
$G$-sequence of $f$ is non-exact at the terms $G_{4}(Y, X;f)$,
$G_{8}^{rel}(Y, X;f)$, and $G_{11}(X)$.

First, $\mathcal{M}_X = \Lambda(x_{4}, x_{11})$, with differential
given on generators by $d(x_4) = 0$, and $d(x_{11}) = x_4^3$, and
$\mathcal{M}_Y = \Lambda(y_8, y_{15}, y_{4}, y_{19})$ with
differential $d(y_8) = 0$, $d(y_{15}) = y_{8}^2$, $d(y_{4}) = 0$,
and $d(y_{19}) = y_{4}^5$. In both models, subscripts denote
degrees. Then the Sullivan model of $f$, which we denote by $\phi
\colon \mathcal{M}_Y \to \mathcal{M}_X$, is given on generators by
$\phi(y_8) = x_{4}^2$, $\phi(y_{15}) = x_{4}x_{11}$, $\phi(y_{4})
= x_{4}$, and $\phi(y_{19}) = x_{4}^2x_{11}$.

For degree reasons, $\Der_i(\mathcal{M}_{X}, \mathcal{M}_{X}; 1) =
0$ unless $i = 3, 4, 7$ or $11$.  Furthermore,
$\Der_{*}(\mathcal{M}_{X}, \mathcal{M}_{X}; 1)$ is spanned by the
derivations $x_{4}^2 \partial x_{11}$, $x_{4}^*$, $x_{4}
\partial x_{11}$, and $x_{11}^*$ of degree $3$, $4$, $7$, and $11$ respectively.
An easy computation reveals that $\delta(x_{4}^*) = -3 x_{4}^2
\partial x_{11}$, but that $x_{4}\partial
x_{11}$ and $x_{11}^*$ are both (non-exact) cocycles. It follows
from \thmref{thm:equivalence of squares} that
\begin{displaymath}
\pi_i\big(\map(X, X; 1)\big)\otimes\Q \cong H_i\big(\Der(\mathcal{M}_{X}, \mathcal{M}_{X}; 1)\big) = \begin{cases} \Q & \text{if $i = 7, 11$} \\
0 & \text{otherwise}  \end{cases}
\end{displaymath}
Further, it is evident that
\begin{displaymath}
\pi_i(X)\otimes\Q \cong H_i\big(\Der(\mathcal{M}_{X}, \Q; \varepsilon)\big) = \begin{cases} \Q & \text{if $i = 4, 11$} \\
0 & \text{otherwise}  \end{cases}
\end{displaymath}
with the non-zero cohomology in degrees $4$ and $11$ generated by
cocycles $x_{4}^*$ and $x_{11}^*$, respectively. Given these
generators, we see that
$$H(\varepsilon_*) \colon H_i\big(\Der(\mathcal{M}_{X},
\mathcal{M}_{X}; 1)\big) \to H_i\big(\Der(\mathcal{M}_{X}, \Q;
\varepsilon)\big),$$
that is, the homomorphism $\omega_{\#}\otimes\Q \colon
\pi_i\big(\map(X, X; 1)\big)\otimes\Q \to \pi_i(X)\otimes\Q$
induced by the evaluation map on rational homotopy groups, is an
isomorphism in degree $11$ and is zero in all other degrees. It
follows from \thmref{thm:equivalence of squares}---see
\corref{cor:rational Gottlieb group}---that $G_i(\mathcal{M}_{X})
= 0$ other than in degree $11$, where we have
$G_{11}(\mathcal{M}_{X}) \cong \Q$. Up to this point, our
observations are both well-known, and also easily obtained by a
number of standard methods.

We now show that the rationalized $G$-sequence is non-exact at the
$G_{11}(\mathcal{M}_{X})$ term.  Recall that this term of the
$G$-sequence, together with its adjacent terms, is obtained from
the diagram
\begin{displaymath}
\xymatrix{\Rel_{12}(\phi^*) \ar[r]^-{P} \ar[d]^-{(\varepsilon_*,
\varepsilon_*)}& \Der_{11}(\mathcal{M}_X, \mathcal{M}_X; 1)
\ar[r]^-{\phi^*}\ar[d]^-{\varepsilon_*} & \Der_{11}(\mathcal{M}_Y,
\mathcal{M}_X; \phi)
\ar[d]^-{\varepsilon_*}\\
\Rel_{12}(\widehat{\phi^*}) \ar[r]^-{\widehat{P}} &
\Der_{11}(\mathcal{M}_X, \Q; \varepsilon)
\ar[r]^-{\widehat{\phi^*}} & \Der_{11}(\mathcal{M}_Y, \Q;
\varepsilon), }
\end{displaymath}
by passing to homology and then considering the image of the top
sequence in  the bottom.  A brute force calculation will display
the result, but we opt to argue at a more general level so as to
indicate some reason for non-exactness.  It is evident that
$H\big(\widehat{\phi^*}\big)\circ H(\varepsilon_*) ([x_{11}^*]) =
0 \in H_{11}\big(\Der(\mathcal{M}_Y, \Q;
\varepsilon)\big)$---indeed, this latter term is zero, since it is
isomorphic to $\pi_{11}(Y)\otimes\Q = 0$. The key point for
non-exactness here, however, is that in the top sequence we have
$H(\phi^*)([x_{11}^*]) \not= 0 \in H_{11}\big(\Der(\mathcal{M}_Y,
\mathcal{M}_X; \phi)\big)$.  In fact, a straightforward check
shows that $\phi^*(x_{11}^*) = x_{4}\partial y_{15} + x_{4}^2
\partial y_{19}$. Since $\Der_{12}(\mathcal{M}_Y, \mathcal{M}_X; \phi)
= 0$, there are no non-zero boundaries in degree $11$ and hence
$H(\phi^*)([x_{11}^*]) \not= 0$.  Consequently, $[x_{11}^*]$
cannot be in the image of $H(P)$ in the top sequence.  Therefore,
since $H(\varepsilon_*)$ is an isomorphism in degree $11$,
$H(\varepsilon_*)([x_{11}^*]) = [x_{11}^*]$ cannot be in the image
of $H(\widehat{P})\circ H(\varepsilon_*, \varepsilon_*)$. It
follows from these facts that $[x_{11}^*] \in
G_{11}(\mathcal{M}_{X})$ is a non-zero element in the kernel of
$H(\widehat{\phi^*})$ and yet is not in the image of
$H(\widehat{P}) \colon G^{rel}_{12}(\mathcal{M}_{Y},
\mathcal{M}_{X}; \phi) \to G_{11}(\mathcal{M}_{X})$.

Next consider the term $G_{4}(\mathcal{M}_{Y}, \mathcal{M}_{X};
\phi)$: Before passing to homology, the relevant diagram is the
following:
\begin{displaymath}
\xymatrix{\Der_{4}(\mathcal{M}_X, \mathcal{M}_X; 1)
\ar[r]^-{\phi^*}\ar[d]^-{\varepsilon_*} & \Der_{4}(\mathcal{M}_Y,
\mathcal{M}_X; \phi) \ar[d]^-{\varepsilon_*} \ar[r]^-{J}&
\Rel_{4}(\phi^*)
 \ar[d]^-{(\varepsilon_*, \varepsilon_*)}
\\
\Der_{4}(\mathcal{M}_X, \Q; \varepsilon)
\ar[r]^-{\widehat{\phi^*}} &  \Der_{4}(\mathcal{M}_Y, \Q;
\varepsilon) \ar[r]^-{\widehat{J}} & \Rel_{4}(\widehat{\phi^*}) ,
}
\end{displaymath}
The derivation $\theta = y_{4}^* + 5 x_{4}x_{11} \partial y_{19}
\in \Der_{4}(\mathcal{M}_Y, \mathcal{M}_X; \phi)$ is a cocycle, as
is easily checked.  Under $H(\varepsilon_*) \colon
H_{4}\big(\Der(\mathcal{M}_Y, \mathcal{M}_X; \phi)\big) \to
H_{4}\big(\Der(\mathcal{M}_Y, \Q; \varepsilon )\big)$, we have
$H(\varepsilon_*)([\theta]) = [y_{4}^*] \not=0$.  Since $[y_{4}^*]
= H(\widehat{\phi^*})([x_{4}^*])$, it follows that
$H(\widehat{J})([y_{4}^*]) = 0$.  As we noted above, however,
$G_{4}(\mathcal{M}_X) = 0$.  Therefore, $[y_{4}^*] \in
G_{4}(\mathcal{M}_{Y}, \mathcal{M}_{X}; \phi)$ is a non-zero
element in the kernel of $H(\widehat{J}) \colon
G_{4}(\mathcal{M}_{Y}, \mathcal{M}_{X}; \phi) \to
G^{rel}_{4}(\mathcal{M}_{Y}, \mathcal{M}_{X}; \phi)$ that is not
in the image of $H(\widehat{\phi^*}) \colon G_{4}(\mathcal{M}_{X})
\to G_{4}(\mathcal{M}_{Y}, \mathcal{M}_{X}; \phi)$.

Finally, consider the term $G^{rel}_{8}(\mathcal{M}_{Y},
\mathcal{M}_{X}; \phi)$: Here, we begin with the following
diagram:
\begin{displaymath}
\xymatrix{\Der_{8}(\mathcal{M}_Y, \mathcal{M}_X; \phi)
\ar[d]^-{\varepsilon_*} \ar[r]^-{J} & \Rel_{8}(\phi^*)
 \ar[d]^-{(\varepsilon_*, \varepsilon_*)} \ar[r]^-{P} & \Der_{7}(\mathcal{M}_X,
\mathcal{M}_X; 1) \ar[d]^-{\varepsilon_*} \\
\Der_{8}(\mathcal{M}_Y, \Q; \varepsilon) \ar[r]^-{\widehat{J}} &
\Rel_{8}(\widehat{\phi^*})  \ar[r]^-{\widehat{P}} &
\Der_{7}(\mathcal{M}_X, \Q; \varepsilon) }
\end{displaymath}
We find that $(-2 x_{4} \partial x_{11}, y_{8}^* + 2 x_{11}
\partial y_{19}) \in \Rel_{8}(\phi^*)$ is a cocycle that has non-zero image in
$H(\varepsilon_*, \varepsilon_*) \colon
H_{8}\big(\Rel(\phi^*)\big) \to
H_{8}\big(\Rel(\widehat{\phi^*})\big)$.  Furthermore, it is
evident that $H(\widehat{P})\circ H(\varepsilon_*,
\varepsilon_*)\big((-2 x_{4} \partial x_{11}, y_{8}^* + 2 x_{11}
\partial y_{19})\big) = H(\widehat{P})([(0, y_{8}^*)]) = 0$---indeed,
$H_{7}\big(\Der(\mathcal{M}_X, \Q; \varepsilon )\big) \cong
\pi_{7}(X)\otimes\Q = 0$.  To see that $[(0, y_{8}^*)]$ is not in
the image of $H(\widehat{J}) \colon G_{8}(\mathcal{M}_{Y},
\mathcal{M}_{X}; \phi) \to G^{rel}_{8}(\mathcal{M}_{Y},
\mathcal{M}_{X}; \phi)$, we will compute $G_{8}(\mathcal{M}_{Y},
\mathcal{M}_{X}; \phi)$ to be zero.

A general derivation $\gamma \in \Der_{8}(\mathcal{M}_{Y},
\mathcal{M}_{X}; \phi)$ can be written as
$$\gamma = \lambda_1 y_{8}^* + \lambda_2 x_{11} \partial y_{19},$$
for coefficients $\lambda_i \in \Q$.  To find the cocycles of this
form, we first observe that $\delta(\gamma)$ must vanish on the
generators $y_{8}$ and $y_{4}$, for degree reasons. On $y_{15}$,
we compute as follows:
$$
\begin{aligned}
\delta(\gamma)(y_{15}) &= (d(\gamma) - (\gamma)d)(y_{15})\\
&= 0 - \gamma(y_{8}^2)\\
&=  - 2 \gamma(y_{8})\phi(y_{8})\\
&=  - 2\lambda_1 x_{4}^2.
\end{aligned}
$$
Thus, if $\gamma$ is a cocycle, then we must have $\lambda_1=0$. A
similar computation shows that $\delta(\gamma)(y_{19}) = \lambda_2
x_{4}^3$, and thus that there are no non-zero cocycles in
$\Der_{8}(\mathcal{M}_{Y}, \mathcal{M}_{X}; \phi)$. In summary, we
have computed that
$$
\pi_{8}\big(\map(X, Y; f)\big)\otimes\Q \cong
H_{8}\big(\Der(\mathcal{M}_{Y}, \mathcal{M}_{X}; \phi)\big) = 0.
$$
This last part of our computation is easily confirmed using the
result of \corref{cor:Grivel}. It follows, of course, that
$G_{8}(\mathcal{M}_{Y}, \mathcal{M}_{X}; \phi) = 0$ and, in
particular, that $[(0, y_{8}^*)]$ is not in the image of
$H(\widehat{J}) \colon G_{8}(\mathcal{M}_{Y}, \mathcal{M}_{X};
\phi) \to G^{rel}_{8}(\mathcal{M}_{Y}, \mathcal{M}_{X}; \phi)$.
\end{example}

\begin{remarks}
The first example of a non-exact $G$-sequence,  given in
\cite{Pak-W}, was in dimension one.  A higher dimensional example
was produced later in \cite{P-S-W}. With the approach illustrated
in the above example, it is straightforward to produce higher
dimensional examples of non-exact rationalized $G$-sequences.
Observe that non-exact ordinary $G$-sequences are produced as a
result, since if a sequence of abelian groups is not exact after
tensoring with $\Q$ then it was not exact to begin with.  Thus,
non-exactness rationally implies non-exactness integrally.
\end{remarks}

Although the $G$-sequence in general is non-exact, there are
certain situations in which it is perfectly well behaved, at least
after rationalization. We now mention one such situation.  We say
that a space $X$ is an \emph{$H_{0}$-space} if its rational
cohomology algebra is a free graded algebra (exterior algebra on
the odd-degree generators tensored with a polynomial algebra on
the even-degree generators). Equivalently, we could require that
$X$ be an $H$-space after rationalization, whence such a space is
also referred to as a \emph{rational $H$-space}.  Recall also the
definition of an $F_{0}$-space from before \corref{cor:Grivel}.

\begin{theorem}\label{thm:exact G-seq}
Let $f \colon X \rightarrow Y$ be any map from an $F_{0}$-space
$X$ to an $H_{0}$-space $Y$ that induces the zero homomorphism on
rational homotopy groups. Then the rationalized $G$-sequence
splits into short exact sequences
$$\xymatrix{0 \ar[r] & G_{n+1}(Y, X;f)\otimes\Q \ar[r]^-{J}
& G^{rel}_{n+1}(Y, X;f)\otimes\Q \ar[r]^-{P} & G_n(X)\otimes\Q
\ar[r] & 0}$$
for each $n \geq 2$.
\end{theorem}

\begin{proof} Our assumption that $f_\#\otimes\Q = 0$ means that the long exact
sequence induced by $f$ on rational homotopy groups splits.
Furthermore, since $Y$ is an $H_0$-space, we have
$G_{n}(Y)\otimes\Q = \pi_n(Y)\otimes\Q$ for each $n$.  It follows
that $G_{n}(Y, X;f)\otimes\Q = \pi_n(Y)\otimes\Q$ for each $n$.
From these general considerations, we have short exact sequences
\begin{equation}\label{eq:short exact general}
\xymatrix{0 \ar[r] & G_{n+1}(Y, X; f )\otimes\Q \ar[r]^-{J} &
\pi_{n+1}(f)\otimes\Q \ar[r]^-{P} & \pi_n(X)\otimes\Q \ar[r] &
0}
\end{equation}
for $n \geq 2$. To sharpen this to the statement of the theorem,
we work within our minimal model framework.

We use some results of Halperin \cite{Hal}, on the rational
homotopy of an $F_{0}$-space.  These state that $X$ has minimal
model $(\mathcal{M}_{X}, d_{X})$ of the form $\mathcal{M}_{X} =
\Lambda(x_{1}, \ldots, x_{n}) \otimes \Lambda(y_1, \ldots, y_{n})$
where $|x_{i}|$ is even, $|y_{j}|$ is odd, $d_{X}(x_{i}) = 0$ and
$d_{X}(y_{j}) \in \Lambda(x_{1}, \ldots, x_{n})$.  Furthermore,
the cohomology is evenly graded, and any cocycle in $I(y_1,
\ldots, y_{n})$, the ideal of $\mathcal{M}_{X}$ generated by the
$y_i$, is exact.  It is well-known that the minimal model of an
$H_0$-space takes the form $\mathcal{M}_{Y} = \Lambda(z_{1},
z_{2}, \ldots)$, with trivial differential. The map $f \colon X
\to Y$ has Sullivan minimal model $\phi \colon \mathcal{M}_Y \to
\mathcal{M}_X$ that is determined by the $\phi(z_i)$.  Since each
$z_i$ is a cocycle, it follows that each $\phi(z_i) \in
\mathcal{M}_X$ is a cocycle. From the results of Halperin
mentioned above, we can write $\phi(z_i) = \chi_i + d(\xi_i)$, for
suitable $\chi_i \in \Lambda(x_{1}, \ldots, x_{n})$ and $\xi_i \in
I(y_1, \ldots, y_{n})$.  The assumption that $f$ induces zero on
rational homotopy groups translates into the further restriction
that each $\chi_i$ is decomposable.

The short exact sequences (\ref{eq:short exact general}),
translated into our derivation setting, correspond to short exact
sequences
$$\xymatrix{0 \to G_{n+1}(\mathcal{M}_{Y}, \mathcal{M}_{X}; \phi) \ar[r]^-{H(\widehat{J})}
& H_{n+1}\big(\Rel(\widehat{\phi^*})\big) \ar[r]^-{H(\widehat{P})}
& H_n\big(\Der(\mathcal{M}_{X}, \Q; \varepsilon)\big) \to 0}$$
We first prove that each map $H(\widehat{P})$ restricts to give a
surjection
$$H(\widehat{P}) \colon G^{rel}_{n+1}(\mathcal{M}_{Y},
\mathcal{M}_{X}; \phi) \to G_{n}(\mathcal{M}_{X}).$$
The Gottlieb elements of $\mathcal{M}_{X}$ are precisely the
$y_{j}^{*}$, dual to the odd-degree generators of
$\mathcal{M}_{X}$.  This can be seen from the description of the
Gottlieb elements given in \corref{cor:rational Gottlieb group},
together with the facts about the minimal model of $X$ recalled
earlier. Now write $\phi(z_i) = \chi_i + d(\xi_i)$ as above and,
for each $j$, define a derivation
$$\theta_j = - \sum_i y_j^*(\xi_i)\partial z_i \in \Der_{|y_j|+1}(\mathcal{M}_{Y}, \mathcal{M}_{X};
\phi).$$
Then $\delta(\theta_j) = - \sum_i d(y_j^*(\xi_i))\partial z_i$,
since the differential in $\mathcal{M}_{Y}$ is trivial.  On the
other hand, we have
$$
\begin{aligned}
\phi^*(y_j^*)(z_i) &= y_i^*\big(\phi(z_i)\big) = y_i^*\big(\chi_i + d(\xi_i)\big) = 0 + y_i^*\big(d(\xi_i)\big)\\
&=  - d\big(y_i^*(\xi_i)\big),
\end{aligned}
$$
with the last step following because $y_j^*$ is a $\delta$-cycle.
Consequently, $(y_j^*, \theta_j) \in \Rel_{|y_j|+1}(\phi^*)$ is a
$\delta_{\phi^*}$-cycle.  Since $H(\widehat{P})\circ
H(\varepsilon_*,\varepsilon_*)(y_j^*, \theta_j) = y_j^*$, it
follows that $H(\widehat{P})$ does restrict to the desired
surjection.

The map $H(\widehat{J})$ is injective on rational homotopy groups,
as we have already observed, and therefore restricts to an
injection in the rationalized $G$-Sequence. So it only remains to
show exactness at the $G^{rel}_{n+1}(Y, X;f)$ terms. Suppose that
$$H(\varepsilon_*, \varepsilon_*)([(\theta, \psi)])
\in \ker\big(H(\widehat{P}) \colon G^{rel}_{n+1}(\mathcal{M}_{Y},
\mathcal{M}_{X}; \phi) \to G_{n}(\mathcal{M}_{X})\big),$$
for some cocycle $(\theta, \psi) \in \Rel_{n+1}(\phi^*)$.  The
fact that $H(\widehat{P})([(\varepsilon_*(\theta) ,
\varepsilon_*(\psi))]) = 0$ implies that $\varepsilon_*(\theta) =
0$.  Thus $(\varepsilon_*, \varepsilon_*)(\theta, \psi) = (0,
\varepsilon_*(\psi))$.  Now define a derivation $\overline{\psi}
\in \Der_{n+1}(\mathcal{M}_{Y}, \mathcal{M}_{X}; \phi)$ by setting
$\overline{\psi} = \psi$ on generators of $\mathcal{M}_{Y}$ of
degree $n+1$ and $\overline{\psi} = 0$ on all other generators of
$\mathcal{M}_{Y}$.  It is easily seen that $\overline{\psi}$ is a
cycle. Indeed, $d_X \overline{\psi} = 0$ on all generators of
$\mathcal{M}_{Y}$, since $\overline{\psi}$ has non-zero image only
in degree zero, and $\overline{\psi}d_Y = 0$ since $d_Y = 0$. Thus
$H(\varepsilon_*)([\overline{\psi}]) \in G_{n+1}(\mathcal{M}_{Y},
\mathcal{M}_{X}; \phi)$ satisfies
$$H(\widehat{J}) \circ H(\varepsilon_*)([\overline{\psi}]) = [(0, \varepsilon \circ \overline{\psi})] =
[(0, \varepsilon \circ \psi)] =
H(\varepsilon_*,\varepsilon_*)([(\theta, \psi)]).$$
That is, $\ker\big(H(\widehat{P})\big) \cap
G^{rel}_{n+1}(\mathcal{M}_{Y}, \mathcal{M}_{X}; \phi)  \subseteq
H(\widehat{J})\big( G_{n+1}(\mathcal{M}_{Y}, \mathcal{M}_{X};
\phi) \big)$ and the rationalized $G$-sequence is exact at each
$G^{rel}_{n+1}(Y, X;f)$ term.
\end{proof}

\begin{remark}
Various conditions are known, under which the $G$-sequence of a
map $f \colon X \to Y$ is exact.  For instance, it is exact when
$f$ is null-homotopic \cite{L-W4}, and when $f$ is a homotopy
monomorphism \cite{P-W}.  The hypotheses of \thmref{thm:exact
G-seq} are well-suited for rational homotopy theory.  Both types
of space are well-known, and it is easy to give examples to which
the theorem applies.   In fact, for fixed $X$ and $Y$, the maps to
which it applies are classified up to rational homotopy by the
decomposable rational cohomology of $X$ that occurs in those
(even) degrees in which $Y$ has a generator of rational
cohomology. We emphasize that the $H_0$-space $Y$ must be allowed
to have polynomial generators in rational cohomology, and hence be
infinite-dimensional, otherwise the theorem reduces to the case in
which the map $f$ is rationally null-homotopic. Furthermore, the
hypothesis that $f$ be zero on rational homotopy groups is
necessary.  For example, the map $f \colon S^4 \to \H P^{\infty}$,
given by inclusion of the bottom cell, has a non-exact
rationalized $G$-sequence, as is easily confirmed by computations
similar to those of \exref{ex:non-exact G-seq}
\end{remark}

Since  the $G$-sequence of a map $f \colon X \rightarrow Y$  is a
boundary sequence, but not usually an  exact sequence, it is
natural to consider its homology. This gives the so-called
$\omega$-homology of $f$ \cite{L-W4}.   In general, one obtains an
$\omega$-homology group at each of the three types of term.  In
the following, we restrict our attention to the $\omega$-homology
that occurs at the Gottlieb group term $G_{*}(X)$, denoted
$H^{a\omega}_{*}(X, Y; f)$ in \cite{L-W4}.  Thus we consider the
sub-quotients of the Gottlieb groups $G_{*}(X)$ defined by
$$ H_{n}^{a\omega}(X,Y;f) = \frac{\ker \{ f_{\#}\colon G_{n}(X) \rightarrow
G_{n}(Y, X;f)   \} }{\im \{P \colon G^{rel}_{n+1}(Y,X;f)
\rightarrow G_{n}(X)  \}}.$$
When $Y$ is an $H_{0}$-space, the rational $\omega$-homology of
$f\colon X\rightarrow Y$ is related to the the negative
derivations on the rational cohomology  of $X$ that are induced by
derivations on the minimal model. To be precise, define a linear
map of degree zero
$$\varphi_{X} \colon H_{*}\big(\Der(\mathcal{M}_{X}, \mathcal{M}_{X};1)\big)
\rightarrow \Der_{*}(H^{*}(X; \Q), H^{*}(X; \Q); 1)$$
by the rule $\varphi_{X}([\theta])([\chi]) = [\theta(\chi)]$, for
$\theta$ a cycle in $\Der_*(\mathcal{M}_{X}, \mathcal{M}_{X};1)$
and $\chi$ a cocycle in $\mathcal{M}_{X}$.  It is straightforward
to check that $\varphi_{X}$ is well-defined. (cf.~
\cite[Prop.1.6]{Griv}.  In fact, $\varphi_{X}$ is a morphism of
graded Lie algebras.)

In the next result, and the example that follows it, we illustrate
that the rationalized $G$-sequence may be exact at all occurrences
of one type of term, while failing to be exact at the other types
of term.  In other words, the rational $\omega$-homology of a map
may be zero at one type of term, yet non-zero at the other types
of term.

\begin{theorem}\label{thm:H^omega=0}
Let $X$ be a finite complex for which the map $\varphi_{X}$
defined above is trivial and let $Y$ be an $H_{0}$-space. Then
$H_{*}^{a\omega}(X,Y;f) \otimes \Q = 0$ for any map $f\colon X
\rightarrow Y$.
\end{theorem}

\begin{proof}    Since $Y$ is an $H_{0}$-space,
its minimal model $\mathcal{M}_{Y} \cong H^*(Y;\Q)$ has trivial
differential. Let $\phi \colon H^*(Y;\Q) \to \mathcal{M}_{X}$
denote the minimal model of $f$. For a derivation $\theta \in
\Der_*(H^*(Y;\Q), \mathcal{M}_{X}; \phi)$, we have $\delta(\theta)
= \pm d_X \theta$.  Using this observation, we obtain a map of
chain complexes
$$\mu \colon \Der_*(H^*(Y;\Q),\mathcal{M}_{X}; \phi) \to \Der_*\big(H^*(Y;\Q),
H^*(X;\Q); H(\phi)\big),$$
defined by $\mu(\theta)(\chi) = [\theta(\chi)]$.  Using the
preceding observation, together with the free-ness of $H^*(Y;\Q)$,
it is straightforward to check that $\mu$ induces an isomorphism
on passing to homology (note that the right-hand term has trivial
differential, and so is its own homology).  Furthermore, the
following diagram commutes:
$$\xymatrix{H_{*}\big(\Der(\mathcal{M}_{X},
\mathcal{M}_{X};1)\big)\ar[rr]^-{H(\phi^*)} \ar[d]_{\varphi_{X}} &
& H_{*}\big(\Der(H^*(Y;\Q), \mathcal{M}_{X};\phi)\big)
\ar[d]^{\mu}_{\cong}
\\
\Der_{*}(H^*(X;\Q), H^*(X;\Q);1) \ar[rr]_{\big(H(\phi)\big)^*}  &
& \Der_{*}\big(H^*(Y;\Q), H^*(X;\Q);H(\phi)\big) }$$
Therefore, the assumption that $\varphi_{X} = 0$ implies that the
top map $H(\phi^*)$ in the above diagram is zero.  A
straightforward diagram chase using the homology ladder that
defines the rationalized $G$-sequence now gives the result.
\end{proof}

\begin{example}
Following \thmref{thm:exact G-seq} we remarked that the cellular
inclusion $S^4 \to \H P^{\infty}$ does not have an exact
rationalized $G$-sequence.  However, it does satisfy the
hypotheses of \thmref{thm:H^omega=0}, since here $X = S^4$ has the
property that all derivations of the cohomology algebra are
trivial.
\end{example}

\begin{remark}
The hypothesis on $X$ in \thmref{thm:H^omega=0}, that
$\varphi_{X}= 0$, deserves some comment.  First, we observe that
the nature of the hypothesis distinguishes structure at the
minimal model level from structure at the cohomology level.  This
is a distinction that is made in rational homotopy for a wide
variety of structures. Next, we observe that this condition is
satisfied for many, if not all, $F_0$-spaces $X$.  Indeed,
\conjref{conj:Halperin}---the long-standing conjecture of Halperin
concerning $F_0$-spaces---is equivalent to the assertion that all
negative-degree derivations on the cohomology algebra of an
$F_0$-space are trivial (see \cite{Me} for details).  Whenever
this conjecture is true---and it has been verified in many
cases---obviously we have $\varphi_{X}= 0$.  Therefore,
\thmref{thm:H^omega=0} can be compared with \thmref{thm:exact
G-seq}, as a result with weaker hypotheses, and correspondingly
weaker conclusion. Finally, we note that the map $\varphi_{X}$
makes an appearance in a completely different context, in the work
of Belegradek and Kapovitch \cite{Be-Ka}.
\end{remark}

Our last set of results relate directly to
\conjref{conj:Halperin}. First, we observe that for an inclusion
of a summand of a product, the $G$-sequence behaves in a
particularly nice way.  Since it is no harder to do so, we state
and prove this result in the integral setting.

\begin{proposition}\label{prop:split short G}
Suppose that $i_1 \colon X \to X \times B$ is the inclusion into
the first summand.  Then the $G$-sequence of $i_1$ is exact, and
furthermore reduces to split short exact sequences
$$
\xymatrix{0 \ar[r] & G_n(X) \ar[r]^-{(i_1)_\#} & G_n(X \times B,
X; i_2) \ar[r]_-{(p_2)_\#} & \pi_n(B) \ar[r] \ar@/_1pc/[l]_-{(i_2)_\#}
 & 0, }
$$
where $p_2 \colon X \times B \to B$ is projection onto the second
summand and the splitting is induced by inclusion into the second
summand  $i_2 \colon B \to X \times B$.
\end{proposition}

\begin{proof}
This follows from results in \cite{L-W1} (see also \cite{Woo}),
but we give a brief argument here. First, let $X
\stackrel{j}{\rightarrow} E \stackrel{p}{\rightarrow} B$ be any
fibre sequence.    Hilton's excision homomorphism for relative
homotopy groups gives an isomorphism $\pi_{*}(j) \cong \pi_{*}(B)$
\cite[Chap.3]{Hilton}. Thus we may view $G^{rel}_{*}(E,X;j)$ as a
subgroup of $\pi_*(B)$ and the $G$-sequence of the fibre inclusion
$j$ as a subsequence the long exact homotopy sequence of the
fibration.

Now apply this remark to the trivial fibration $X
\stackrel{i_1}{\rightarrow} X \times B \stackrel{p_2}{\rightarrow}
B$.  The inclusion $i_2 \colon B \to X \times B$ induces a
splitting of the long exact homotopy sequence of this trivial
fibration in the usual way.  The result now follows from the
observation that $(i_2)_\#\big(\pi_n(B)\big) \subseteq G_n(X
\times B, X;p_2)$, as is easily established from the definitions.
\end{proof}

Of course, this result and its proof can be rationalized, and it
is in the rational setting that we will use it.
\conjref{conj:Halperin} concerns fibrations $X \to E \to B$ with
fibre an $F_0$-space and arbitrary base.  However, it is
well-known how to reduce the conjecture to consideration of such
fibrations with base an odd-dimensional sphere \cite{Me}.
Furthermore, in \cite{Lu} it is pointed out that, for such
fibrations with base an odd-dimensional sphere, Halperin's
conjecture actually asserts that the fibration should be trivial.

From these remarks, we see that a necessary condition for
\conjref{conj:Halperin} to be true is that any fibration $X \to E
\to S^{2r+1}$ with $X$ an $F_0$-space must have a fibre inclusion
whose $G$-sequence reduces to the split short exact sequences
corresponding via \propref{prop:split short G} to the inclusion
$i_1 \colon X \to X \times S^{2r+1}$.  Perhaps surprisingly, the
converse is true.

\begin{theorem}\label{thm:Halperin G-trivial}
Let $X \stackrel{j}{\rightarrow} E \stackrel{p}{\rightarrow}
S^{2r+1}$ be any fibration with $X$ an $F_0$-space.  The following
are equivalent:
\begin{enumerate}
\item The fibration is rationally TNCZ, that is, $j^*\colon H^*(E;\Q) \to H^*(X;\Q)$ is surjective.
\item The rationalized $G$-sequence of the fibre inclusion reduces to split
short exact sequences
$$
\xymatrix{0 \ar[r] & G_n(X) \otimes \Q \ar[r]^-{(j)_\#} & G_n(E, X;
j) \otimes \Q
\ar[r]_-{(p)_\#} & \pi_n(S^{2r+1}) \otimes \Q \ar[r]
\ar@/_1pc/[l]_-{(i_2)_\#}
 & 0}
$$
for each $n \geq 2.$
\end{enumerate}
\end{theorem}

\begin{proof}
The implication (1) $\Rightarrow$ (2) follows by the remarks
preceding the enunciation.  We prove (2) $\Rightarrow$ (1).

Suppose the fibration $X \to E \to S^{2r+1}$ has minimal model
$$\xymatrix{\Lambda(u) \ar[r]^-{i} & \Lambda(u)\otimes\Lambda V, D
\ar[r]^-{\pi} & (\Lambda V, d), }
$$
where $i$ denotes the inclusion $i(u) = u\otimes1$ and $\pi$ is
the projection.  The hypothesis that $p_\# \colon G_{2r+1}(E,
X; j) \otimes \Q \to \pi_{2r+1}(S^{2r+1}) \otimes \Q$ is onto---included in (2), when
translated into our derivation setting, gives the existence of a
$\pi$-derivation $\psi \in \Der_{2r+1}(\Lambda(u)\otimes\Lambda V,
\Lambda V; \pi)$ that is a cocycle, and that satisfies $\psi(u) =
1$.  Using this $\psi$, define a linear map $\Phi \colon
\Lambda(u)\otimes\Lambda V \to \Lambda(u)\otimes\Lambda V$ by
setting $\Phi(a + ub) = a + u b + u \psi(a)$ for a typical element
$a + ub \in \Lambda (u)\otimes\Lambda V$. We will check that
$\Phi$ is actually a DG algebra isomorphism
$(\Lambda(u)\otimes\Lambda V, D) \to (\Lambda(u)\otimes\Lambda V,
d)$.   First,  $\Phi$ is an algebra map. This follows from the
fact that $\psi$ is a derivation. Suppose given two elements $a +
ub, a' + ub' \in \Lambda(u)\otimes\Lambda V$.  Then we have
$$
\begin{aligned}
\Phi\big((a + ub)(a' + ub')\big) &= \Phi\big(aa' + u( (-1)^{|a|}ab' + ba')\big) \\
&= aa' + u( (-1)^{|a|}ab' + ba') + u \psi(aa').
\end{aligned}
$$
On the other hand, we have
$$
\begin{aligned}
\Phi(a + ub) \Phi(a' + ub') &= \big(a + u b + u \psi(a)\big)
\big(a' + u b' + u \psi(a')\big)\\
&=  aa' + u\big( (-1)^{|a|}ab' + ba' + (-1)^{|a|} a \psi(a') +
\psi(a)a'.
\end{aligned}
$$
These agree, since $\psi(aa') = \psi(a)a' + (-1)^{|a|} a
\psi(a')$. Next, we check that  $\Phi$ commutes with
differentials, that is, that $\Phi D = d \Phi$.  This will follow
from the fact that $\psi$ is a cocycle.  First observe that
$d\Phi(u) = 0 = \Phi D(u)$.  It is thus sufficient to check that
$\Phi D(\chi) = d \Phi(\chi)$ for a typical element $\chi \in
\Lambda V$.  To this end, write $D(\chi) = d(\chi) +
u\theta(\chi)$.  Now we calculate
$$
\Phi D \chi = \Phi\big( d(\chi) + u\theta(\chi)\big) = d(\chi) +
u\theta(\chi) + u \psi(d\chi)
$$
and
$$
d \Phi \chi = d \big( \chi + u\psi(\chi)\big) = d(\chi) - u d
\psi(\chi).
$$
These agree since $\psi$ is a cocycle, whence $-d\psi(\chi) = \psi
D(\chi) = \psi\big(d(\chi) + u\theta(\chi)\big)$. It is evident
that $\Phi$ is an isomorphism, since we have $\Phi(u) = u$ and
$\pi\circ\Phi = \pi$. Therefore, $\Phi \colon
(\Lambda(u)\otimes\Lambda V, D) \to (\Lambda(u)\otimes\Lambda V,
d)$ is a DG isomorphism and the fibration is rationally trivial.
\end{proof}

This leads to the following equivalent phrasing of
\conjref{conj:Halperin}:

\begin{corollary}\label{cor:G-sequence Halperin} Let $X$ be an $F_{0}$-space.  Then $X$ satisfies \conjref{conj:Halperin}
 if and only if the $G$-sequence of the fibre
inclusion in every fibration of the form $X \rightarrow E
\rightarrow S^{2n+1}$ decomposes into split short exact sequences
as in (2) of \thmref{thm:Halperin G-trivial}.
\end{corollary}

\begin{proof}
In \cite{Me}, Meier showed that  Halperin's conjecture for $X$ is
equivalent to the collapsing of the rational Serre spectral
sequence for all fibrations with fibre $X$ and base an odd sphere.
By \cite[Th.2.3]{Lu} this latter condition is equivalent to the
rational homotopy triviality of the fibration.  The result now
follows from \thmref{thm:Halperin G-trivial}.
\end{proof}

\begin{remark}
It is possible to develop the ideas leading to
\thmref{thm:Halperin G-trivial} substantially beyond the
application given above.   Namely, one can use the splitting of
the $G$-sequence of the fibre inclusion as a \emph{measure of how
close the fibration is to being trivial}.  Now in the above
results, we see that for fibrations of the form $X \rightarrow E
\rightarrow S^{2n+1}$, with $X$ an $F_0$-space, this notion
coincides with the notion of the fibration being TNCZ.
Furthermore, in this closely restricted setting, both coincide
with the fibration actually being trivial. In general, however,
this $G$-sequence point of view gives a new way of measuring how
close a fibration is to being trivial. We intend to develop these
ideas in a subsequent paper.
\end{remark}

\begin{appendix}

\section{}

In this appendix, we give careful proofs of the results from DG
algebra homotopy theory that are used to establish
\thmref{thm:equivalence of squares}. Since it is a technical
appendix, we rely on a greater degree of familiarity with
techniques from rational homotopy theory.  We use the notion of
pullback in the DG algebra setting. By this, we mean the
following. Suppose given DG algebra maps $f\colon A \to C$ and
$g\colon B \to C$. Then we form the \emph{(DG algebra) pullback}
(or \emph{fibre product}, as it is called in \cite{F-H-T}) as
$A\oplus_C B = \{ (x, y) \in A\oplus B \mid  f(x) = g(y)\}$. Here
$A \oplus B$ denotes the direct sum of DG algebras. Together with
the projections, the pullback forms the following (strictly)
commutative square of DG algebra maps:
\begin{equation}\label{eq:pullback}
\xymatrix{A\oplus_C B \ar[r]^-{p_1}\ar[d]_{p_2} & A \ar[d]^{f}\\
B\ar[r]_{g} & C}
\end{equation}
This square possesses the usual universal property of pullbacks.
Namely, suppose given DG algebra maps $\alpha \colon Z \to A$ and
$\beta \colon Z \to B$ that satisfy  $f\circ\alpha = g
\circ\beta$.  Then there exists a DG algebra map $\phi =
(\alpha,\beta) \colon Z \to A\oplus_C B$, which is the unique DG
algebra map for which $p_1\circ\phi = \alpha$ and $p_2\circ\phi =
\beta$. We emphasize that throughout this appendix we distinguish
carefully between diagrams that are strictly commutative and ones
that are commutative only up to DG homotopy.  Indeed, it is
precisely this distinction that calls for the proofs of this
appendix.

The following basic property of the pullback is readily gleaned
from the discussion in \cite[Sec.13(a)]{F-H-T}:

\begin{lemma}\label{lem:basic pullback}
Suppose that either $f$ or $g$ is surjective in the pullback
diagram (\ref{eq:pullback}).  If $f$ is a quasi-isomorphism, then
$p_2$ is a quasi-isomorphism.
\end{lemma}

We also use the so-called \emph{surjective trick}, described in
\cite[Sec.12(b)]{F-H-T}.  Given a DG algebra map $\eta \colon B
\to A$, this manoeuvre results in a diagram
$$\xymatrix{B \ar[d]_{\eta} \ar@<0.75ex>[r]^-{\lambda} & B\otimes E(A)
\ar@<0.75ex>[l]^-{1\cdot\varepsilon} \ar@{>>}[ld]^{\gamma}\\
A}$$
in which $\gamma$ is a surjection, and both $1\cdot\varepsilon$
and $\lambda$ are quasi-isomorphisms. Some parts of the diagram
commute, thus $(1\cdot\varepsilon)\circ\lambda = 1$ and
$\gamma\circ\lambda = \eta$.   Other compositions result in
commutativity only up to DG homotopy.  Recall that the notion of
DG homotopy is only defined for DG algebra maps from a minimal
model. Given any map $\phi \colon \mathcal{M} \to B\otimes E(A)$
from a minimal model into $B\otimes E(A)$, we have $\phi \sim
\lambda\circ(1\cdot\varepsilon)\circ\phi$, where $\sim$ denotes DG
homotopy of maps from a minimal model. In particular, we thus have
$\eta\circ(1\cdot\varepsilon)\circ\phi \sim \gamma\circ\phi$.

Now suppose given a map $f \colon X \to Y$.  We choose and fix a
minimal model $\mathcal{M}_f \colon \mathcal{M}_Y \to
\mathcal{M}_X$ for $f$ as follows (\cite[Sec.12(c)]{F-H-T}):  Let
$A^*(f) \colon A^*(Y) \to A^*(X)$ denote the map induced by $f$ on
polynomial differential forms. Let $\eta_X \colon \mathcal{M}_X$
$\to A^*(X)$ and $\eta_Y \colon \mathcal{M}_Y \to A^*(Y)$ denote
minimal models for $X$ and $Y$. As in \cite[Sec.12(b)]{F-H-T}, we
convert $\eta_X$ into a surjection $\gamma_X \colon
\mathcal{M}_X\otimes E(A^*(X)) \to A^*(X)$ and lift
$A^*(f)\circ\eta_Y$ through the surjective quasi-isomorphism
$\gamma_X$, using \cite[Lem.12.4]{F-H-T}, to obtain $\phi_f \colon
\mathcal{M}_Y \to \mathcal{M}_X\otimes E(A^*(X))$. Now set
$\mathcal{M}_f = \gamma_X\circ \phi_f$.  All this is summarized in
the following diagram.
\begin{displaymath}
\xymatrix{ & & & \mathcal{M}_X\otimes E(A^*(X))
\ar@{>>}@/^1pc/[ddl]^-{\gamma_X}_(0.6){\simeq}
\ar@<1ex>[dl]^(0.6){\beta}\ar@<1ex>[dl]^(0.4){\simeq} \\
\mathcal{M}_Y \ar@{.>}[rr]_{\mathcal{M}_f}
\ar@/^1pc/[rrru]^-{\phi_{f}}\ar[d]_{\eta_Y}^{\simeq} & &
\mathcal{M}_X \ar[d]_{\eta_X}^{\simeq} \ar[ur]^-{\alpha}\ar[ur]^(0.3){\simeq}\\
 A^*(Y) \ar[rr]_{A^*(f)} & & A^*(X)}
\end{displaymath}
In this and subsequent diagrams, we indicate that a map is a
quasi-isomorphism with the symbol $\simeq$. By construction, we
have $\gamma_X\circ\phi_f = A^*(f)\circ\eta_Y$,
$\gamma_X\circ\alpha = \eta_X$, and $\beta\circ\alpha = 1$.
Remaining parts of the diagram only commute up to DG homotopy,
however, thus we have $\eta_X \circ\beta \sim \gamma_X$,
$\eta_X\circ\mathcal{M}_f \sim A^*(f)\circ\eta_Y$, and so-on.

Now let $\alpha \colon S^n \to \map(X,Y;f)$ be a representative of
a homotopy class in $\pi_n\big(\map(X,Y;f)\big)$.  Let $F \colon
S^n \times X \to Y$ be an affiliated map for $\alpha$, that is,
$F(s, x) = \alpha(s)(x)$. Since $\alpha$ is a based map, we have
$F\circ i = f \colon X \to Y$, where $i\colon X \to S^n \times X$
denotes (based) inclusion into the second summand $i(x) = (*, x)$.

In the following result, we justify that the Sullivan minimal
model of any affiliated map, and a DG homotopy between two such,
have the restricted form that we require of them for the
definition and well defined-ness of $\Phi_f$.

\begin{proposition}\label{prop:restricted DG homotopy}
Suppose given maps $F, G \colon S^n \times X \to Y$ and $H \colon
S^n \times X \times I \to Y$ a homotopy from $F$ to $G$ that is
stationary on $X\times I$. Suppose that $H(*, x, t) = f(x)$ for $f
\colon X \to Y$ and let $\mathcal{M}_f \colon \mathcal{M}_Y \to
\mathcal{M}_X$ be a fixed choice of minimal model for $f$.   There
is a DG homotopy $\mathcal{H} \colon \mathcal{M}_Y \to
\mathcal{M}_{S^n}\otimes \mathcal{M}_X\otimes \Lambda(t, dt)$ from
$\mathcal{M}_F$ to $\mathcal{M}_G$, minimal models for $F$ and
$G$, of the form
\begin{displaymath}
\mathcal{H}(\chi) = 1\otimes \mathcal{M}_f(\chi)\otimes1 +
\mathrm{terms\  in\ } \big(\mathcal{M}_{S^n}\big)^{+}\otimes
\mathcal{M}_X\otimes \Lambda(t, dt).
\end{displaymath}
In particular, any map $F\colon S^n \times X \to Y$ that satisfies
$F\circ i = f$ has a minimal model $\mathcal{M}_F \colon
\mathcal{M}_Y \to \mathcal{M}_{S^n}\otimes \mathcal{M}_X$ of the
form
\begin{displaymath}
\mathcal{M}_F(\chi) = 1\otimes \mathcal{M}_f(\chi) +
\mathrm{terms\  in\ } \big(\mathcal{M}_{S^n}\big)^{+}\otimes
\mathcal{M}_X.
\end{displaymath}
\end{proposition}

\begin{proof}
Construct the following pullback:
\begin{displaymath}
\xymatrix{P \ar[r]^-{p_1}\ar[d]_{p_2} & \mathcal{M}_X \otimes
\Lambda(t, dt) \otimes E(A^*(X\times I))
 \ar@{>>}[d]^{\gamma}_{\simeq}\\
A^*(S^n\times X\times I)\ar@{>>}[r]_-{A^*(i)} & A^*(X\times I)}
\end{displaymath}
Here, $i \colon X\times I \to S^n \times X\times I$ denotes the
inclusion $i(x,t) = (*, x, t)$ and $\gamma$ denotes the surjective
quasi-isomorphism obtained by converting the quasi-isomorphism
$\mathcal{M}_X\otimes \Lambda(t, dt) \to A^*(X\times I)$ to a
surjection. Since $i$ is an inclusion, the induced map $A^*(i)$ is
a surjection. From \lemref{lem:basic pullback}, we have that $p_2$
is a quasi-isomorphism.  Now consider the following commutative
diagram:
\begin{displaymath}
\xymatrix@C=-20pt{ & \mathcal{M}_{S^n}\otimes\mathcal{M}_X\otimes
\Lambda(t, dt)\otimes E(A^*(S^n\times X))
\ar@{.>}[dd]|(0.68)\hole^(0.3){\phi}
\ar@/^1pc/@{>>}[rdd]^(0.7){\;\varepsilon\cdot\big(1\otimes1\otimes
E(A^*(i))\big)}
\ar@/_2pc/@{>>}[ddd]|(0.46)\hole^(0.3){\gamma'} \\
 & &  \\
\mathcal{M}_Y \ar@{.>}[r]|(0.74)\hole_{\psi}
\ar@{.>}[ruu]^{\phi_H} \ar@/_/[dr]_(0.4){A^*(H)\circ\eta_Y}
\ar@/^2pc/[rr]^(0.55){\big(J\otimes E(A^*(i))\big)\circ\phi_f} & P
\ar[r]^-{p_1}\ar[d]^{p_2} & \mathcal{M}_X \otimes \Lambda(t, dt)
\otimes E(A^*(X\times I))
 \ar@{>>}[d]^{\gamma}\\
 & A^*(S^n\times X\times I)\ar[r]_-{A^*(i)} & A^*(X\times I),}
\end{displaymath}
Here, $\gamma'$ is the surjective quasi-isomorphism obtained by
converting the quasi-isomorphism
$\mathcal{M}_{S^n}\otimes\mathcal{M}_X \otimes \Lambda(t, dt)\to
A^*(S^n\times X\times I)$ to a surjection, $\phi_f \colon
\mathcal{M}_Y \to \mathcal{M}_X\otimes E(A^*(X))$ denotes the lift
used above to obtain our fixed choice of minimal model for $f$,
and $J \colon \mathcal{M}_X \to \mathcal{M}_X\otimes \Lambda(t,
dt)$ denotes the map $J(x) = x\otimes1$.  From the pullback, we
obtain the maps $\phi$ and $\psi$ indicated. Since both $p_2$ and
$\gamma'$ are quasi-isomorphisms, it follows that $\phi$ is a
quasi-isomorphism. We claim that $\phi$ is also surjective. For
suppose $(a,b) \in P$, so that $\gamma(a) = A^*(i)(b)$. Since
$\varepsilon\cdot\big(1\otimes1\otimes E(A^*(i))\big)$ is
surjective, we can pick $x \in
\mathcal{M}_{S^n}\otimes\mathcal{M}_X\otimes \Lambda(t, dt)\otimes
E(A^*(S^n\times X\times I))$ with $\phi(x) = (a, b')$, and
$A^*(i)(b') = \gamma(a) = A^*(i)(b)$.  Thus $b-b'\in
\ker\,A^*(i)$.  So now pick $y \in E(A^*(S^n\times X\times I))$
with $\gamma_{S^n}(1\otimes1\otimes1\otimes y) = b - b'$---recall
that $E(A^*(S^n\times X\times I))$ is freely generated by the
vector space $A^*(S^n\times X\times I)$ and its suspension.  Then
$\phi(x + 1\otimes1\otimes1\otimes y) = (a,b)$; so $\phi$ is
indeed surjective. Now lift $\psi$ through the surjective
quasi-isomorphism $\phi$, to obtain a map $\phi_H \colon
\mathcal{M}_{Y} \to \mathcal{M}_{S^n}\otimes\mathcal{M}_X\otimes
\Lambda(t, dt)\otimes E(A^*(S^n\times X\times I))$ with
$\phi\circ\phi_H = \psi$.  As usual for the minimal model of a
map, we now set $\mathcal{H} = \beta\circ\phi_H \colon
\mathcal{M}_{Y} \to \mathcal{M}_{S^n}\otimes\mathcal{M}_X\otimes
\Lambda(t, dt)$.

We now check that this DG homotopy starts and ends at minimal
models for $F$ and $G$ respectively.  Let $j_0 \colon S^n\times X
\to S^n \times X \times I$ denote inclusion at $t=0$, and
$\varepsilon_0 \colon \Lambda(t, dt) \to \Q$ the augmentation
given by $\varepsilon_0(t) = 0$ and $\varepsilon_0(dt) = 0$.  We
must check that $(1\otimes1)\cdot\varepsilon_0 \colon
\mathcal{M}_{S^n}\otimes\mathcal{M}_X\otimes \Lambda(t, dt) \to
\mathcal{M}_{S^n}\otimes\mathcal{M}_X$ is a minimal model for $F$.
For this,  consider the following diagram:
\begin{displaymath}
\xymatrix@C=-12pt{ & \mathcal{M}_{S^n}\otimes\mathcal{M}_X\otimes
\Lambda(t, dt)\otimes E\big(A^*(S^n\times X \times
I)\big)\ar@{>>}@/^5pc/[ddd]|(0.67)\hole ^(0.5){\gamma'}
\ar@<1ex>[dd]^{\beta'}
\ar[rd]^(0.7){\big((1\otimes1)\cdot\varepsilon_0\big)\otimes
E(A^*(j_0))\big)} \\
&  & \mathcal{M}_{S^n}\otimes\mathcal{M}_X\otimes
E\big(A^*(S^n\times X)\big) \ar@{>>}@/^3pc/[dd]^-{\gamma''}
\ar@<1ex>[d]^(0.6){\beta''} \\
\mathcal{M}_Y \ar[r]_-{\mathcal{H}}
\ar@/^1pc/[ruu]^-{\phi_{H}}\ar[d]_{\eta_Y}  &
\mathcal{M}_{S^n}\otimes\mathcal{M}_X\otimes \Lambda(t, dt)
\ar[d]_{\eta'}
\ar[uu]^-{\alpha'}\ar[r]_{(1\otimes1)\cdot\varepsilon_0}
& \mathcal{M}_{S^n}\otimes\mathcal{M}_X \ar[u]^-{\alpha''} \ar[d]_{\eta''}\\
 A^*(Y) \ar[r]_-{A^*(H)} &  A^*(S^n\times X\times I)\ar[r]_{A^*(j_0)} &  A^*(S^n\times X)}
\end{displaymath}
We have
\begin{displaymath}
\begin{aligned}
\eta''\circ (1\otimes1)\cdot\varepsilon_0\circ \mathcal{H} & =
A^*(j_0)\circ\eta'\circ \beta'\circ\phi_H\\
&\sim A^*(j_0)\circ\gamma'\circ\phi_H = A^*(j_0)\circ
A^*(H)\circ\eta_Y \\
&= A^*(F)\circ\eta_Y.
\end{aligned}
\end{displaymath}
It follows that $(1\otimes1)\cdot\varepsilon_0\circ \mathcal{H}$
is a minimal model for $F$.  A similar argument shows that
$(1\otimes1)\cdot\varepsilon_1\circ \mathcal{H}$ is a minimal
model for $G$, where $\varepsilon_1 \colon \Lambda(t, dt) \to \Q$
is the augmentation given by $\varepsilon_1(t) = 1$ and
$\varepsilon_1(dt) = 0$.
\end{proof}

Finally, we show that $\Phi'_f \colon \pi_{n}(\map(X,Y;f)) \to
H_{n}\big(\Der(\mathcal{M}_{Y},\mathcal{M}_{X};\mathcal{M}_{f})\big)$,
as defined in the proof of \thmref{thm:equivalence of squares}, is
a homomorphism.  The key here is to translate addition in
$\pi_{n}(\map(X,Y;f))$ into addition of homotopy classes of
affiliated maps.   Denote by $[S^n \times X, Y]_f$ the set of
homotopy classes of affiliated maps $F$ that restrict to $f\colon
X \to Y$. Now suppose $\alpha, \beta \in \pi_{n}(\map(X,Y;f))$
correspond to $F, G \in [S^n \times X, Y]_f$.  Then the sum
$\alpha+ \beta \in \pi_{n}(\map(X,Y;f))$ corresponds to $F+G \in
[S^n \times X, Y]_f$, described as follows.  First consider the
commutative diagram
\begin{displaymath}
\xymatrix{ & X
 \ar[rr]^-{i_2} \ar[dl]^-{i_2}
\ar'[d]^-{f}[dd] & & S^n\times X \ar[dl]^-{j_2}
\ar[dd]^-{G} \\
S^n\times X \ar[rr]_(0.4){j_1}
 \ar[dd]_{F}& &
P  \ar[dd]^(0.3){(F\mid G)}\\
 &  Y \ar[dl]_{1_Y}
\ar'[r]^-{1_Y}[rr]\ar[dl]^{1_Y} & & Y \ar[dl]^{1_Y} \\
Y  \ar[rr]_{1_Y} & & Y }
\end{displaymath}
in which the top and bottom squares are pushouts, and the front
vertical map is induced as a map of pushouts from the back three.
The pushout $P$ can readily be identified with $(S^n \vee
S^n)\times X$. Let $\sigma \colon S^n \to S^n \vee S^n$ denote the
usual pinching comultiplication.  Then we obtain the composition
$(F\mid G)\circ(\sigma\times1) \colon S^n\times X \to Y$, which
gives $F+G$.

Now we translate this addition in $[S^n \times X, Y]_f$ into
minimal model terms.  Recall the notation established in the proof
of \thmref{thm:equivalence of squares}.

\begin{proposition}\label{prop:Phi homomorphism}
Suppose $F$ and $G$ have minimal models such that
$(\psi\otimes1)\circ\mathcal{M}_F(\chi) =
1\otimes\mathcal{M}_f(\chi) + s_n \otimes \theta_F(\chi)$ and
$(\psi\otimes1)\circ\mathcal{M}_G(\chi) =
1\otimes\mathcal{M}_f(\chi) + s_n \otimes \theta_G(\chi)$.  Then
the map $F+G = (F\mid G)\circ(\sigma\times1)$ has minimal model
such that $(\psi\otimes1)\circ\mathcal{M}_{F+G}(\chi) =
1\otimes\mathcal{M}_f(\chi) + s_n \otimes
(\theta_F+\theta_G)(\chi)$.
\end{proposition}

\begin{proof}
First consider the map $(F\mid G) \colon (S^n \vee S^n)\times X
\to Y$.  As discussed in \cite[Sec.13(b)]{F-H-T}, a DG algebra
model for the pushout $P = (S^n \vee S^n)\times X$ can be obtained
as a pullback of minimal models. Specifically, by forming the
pullback of minimal models
\begin{displaymath}
\xymatrix{Q \ar[r]^-{p_1}\ar[d]_{p_2} &
H^*(S^n;\Q)\otimes\mathcal{M}_X
 \ar@{>>}[d]^{\varepsilon\cdot1}\\
H^*(S^n;\Q)\otimes\mathcal{M}_X\ar@{>>}[r]_-{\varepsilon\cdot1} &
\mathcal{M}_X,}
\end{displaymath}
we obtain a DG algebra $Q = H^*(S^n;\Q)\otimes\mathcal{M}_X
\oplus_{\mathcal{M}_X}H^*(S^n;\Q)\otimes\mathcal{M}_X$
quasi-isomorphic to the minimal model of $P$.  Further, the same
discussion can readily be extended to show that in the above
situation, a minimal model for the induced map $(F\mid G)$ can be
obtained using the map induced on the pullback of minimal models.
Now this induced map is the map
\begin{displaymath}
\xi \colon \mathcal{M}_Y \to H^*(S^n;\Q)\otimes\mathcal{M}_X
\oplus_{\mathcal{M}_X}H^*(S^n;\Q)\otimes\mathcal{M}_X
\end{displaymath}
given by $\xi(\chi) = \big(1\otimes\mathcal{M}_f(\chi) + s_n
\otimes \theta_F(\chi), 1\otimes\mathcal{M}_f(\chi) + s_n \otimes
\theta_G(\chi)\big)$.  Under the identification of the pullback
$Q$ with $\big(H^*(S^n;\Q)\oplus
H^*(S^n;\Q)\big)\otimes\mathcal{M}_X$, the composition $(F\mid
G)\circ(\sigma\times1)$ evidently has minimal model of the
described form.
\end{proof}

\end{appendix}

\providecommand{\bysame}{\leavevmode\hbox
to3em{\hrulefill}\thinspace}



\end{document}